\documentclass[twoside]{article}
\usepackage{amsmath,amssymb,amscd,amsthm,verbatim,alltt,amsfonts,array}
\usepackage{mathrsfs}
\usepackage[english]{babel}
\usepackage{latexsym}
\usepackage{amssymb}
\usepackage{euscript}
\usepackage{graphicx}
\usepackage[colorlinks,plainpages]{hyperref}

\usepackage{caption}

\newcommand{\salt}{\vspace*{2.5mm}}     

\def\boxit#1{\vbox{\hrule\hbox{\vrule\kern2pt
 \vbox{\kern4pt#1\kern2pt}\kern2pt\vrule}\hrule}}
\def\qed{\hfill
 $\hskip0.3cm
 \boxit{\hsize 2pt \vsize 10pt}$\bigskip\noindent}

\def\C{{\mathchoice {\setbox0=\hbox{$\displaystyle\rm C$}\hbox{\hbox
to0pt{\kern0.4\wd0\vrule height0.9\ht0\hss}\box0}}
{\setbox0=\hbox{$\textstyle\rm C$}\hbox{\hbox
to0pt{\kern0.4\wd0\vrule height0.9\ht0\hss}\box0}}
{\setbox0=\hbox{$\scriptstyle\rm C$}\hbox{\hbox
to0pt{\kern0.4\wd0\vrule height0.9\ht0\hss}\box0}}
{\setbox0=\hbox{$\scriptscriptstyle\rm C$}\hbox{\hbox
to0pt{\kern0.4\wd0\vrule height0.9\ht0\hss}\box0}}}}
\def\Q{{\mathchoice {\setbox0=\hbox{$\displaystyle\rm
Q$}\hbox{\raise
0.15\ht0\hbox to0pt{\kern0.4\wd0\vrule height0.8\ht0\hss}\box0}}
{\setbox0=\hbox{$\textstyle\rm Q$}\hbox{\raise
0.15\ht0\hbox to0pt{\kern0.4\wd0\vrule height0.8\ht0\hss}\box0}}
{\setbox0=\hbox{$\scriptstyle\rm Q$}\hbox{\raise
0.15\ht0\hbox to0pt{\kern0.4\wd0\vrule height0.7\ht0\hss}\box0}}
{\setbox0=\hbox{$\scriptscriptstyle\rm Q$}\hbox{\raise
0.15\ht0\hbox to0pt{\kern0.4\wd0\vrule height0.7\ht0\hss}\box0}}}}
\def\T{{\mathchoice {\setbox0=\hbox{$\displaystyle\rm
T$}\hbox{\hbox to0pt{\kern0.3\wd0\vrule height0.9\ht0\hss}\box0}}
{\setbox0=\hbox{$\textstyle\rm T$}\hbox{\hbox
to0pt{\kern0.3\wd0\vrule height0.9\ht0\hss}\box0}}
{\setbox0=\hbox{$\scriptstyle\rm T$}\hbox{\hbox
to0pt{\kern0.3\wd0\vrule height0.9\ht0\hss}\box0}}
{\setbox0=\hbox{$\scriptscriptstyle\rm T$}\hbox{\hbox
to0pt{\kern0.3\wd0\vrule height0.9\ht0\hss}\box0}}}}
\def\Z{{\mathchoice {\hbox{$\sf\textstyle Z\kern-0.4em Z$}}
{\hbox{$\sf\textstyle Z\kern-0.4em Z$}}
{\hbox{$\sf\scriptstyle Z\kern-0.3em Z$}}
{\hbox{$\sf\scriptscriptstyle Z\kern-0.2em Z$}}}}

\newcommand{\eqnoone}  
   {}
\newcommand{\eqnotwo}  
   {}

\newcounter{alf}

\newcommand{\adresa}[1]{\par\vspace*{-11pt}
                        \begin{flushright}
                        {\small
                        #1}
                        \end{flushright}
                        }

\def\N{\Bbb N}
\def\Z{\Bbb Z}
\def\R{\Bbb R}
\def\C{\Bbb C}

\newtheorem{Theorem}{Theorem}
\theoremstyle{plain}

\newtheorem{Corollary}{Corollary}

\newtheorem{Definition}{Definition}
\newtheorem{Example}{Example}

\newtheorem{Lemma}{Lemma}

\newtheorem{Proposition}{Proposition}
\newtheorem{Remark}{Remark}

\numberwithin{equation}{section}

\DeclareMathOperator{\ord}{ord}

\DeclareMathOperator{\Spec}{Spec}
\DeclareMathOperator{\Max}{Max}

\newcommand{\K}{{\mathbb K}} 
\newcommand{\fm}{\mathfrak{m}}
\newcommand{\eps}{\varepsilon}

\newcolumntype{L}[1]{>{\raggedright\arraybackslash}p{#1}}
\newcommand{\stack}[1]{\begin{tabular}[t]{@{}l@{}}#1\end{tabular}}

\begin{document}

\begin{center} {\bf Unimodal Singularities of Parametrized Space Curves}\\
          {by}\\
{\sc Gert-Martin Greuel, Khawar Mehmood and Gerhard Pfister}
\end{center}

\begin{abstract}
The aim of the article is to classify unimodal parametrized space curve singularities defined over an algebraically closed field of characteristic $0$. We provide explicit normal forms of parametrizations up to $\mathcal A$-equivalence. The parametrizations are given by convergent power series over a 
 real valued field, which includes parametrizations given by formal power series and by, for example, convergent power series over the complex numbers or the field of Puisieux series. 
The list of normal forms is surprisingly simple.

\end{abstract}

\begin{quotation}

\end{quotation}

\thispagestyle{empty}

\section{Introduction and results}\label{Sec1}
The study and classification of singularities have a long history. Important contributions go back to V.I. Arnold \cite{A1}, 
\cite{A2} for hypersurface singularities and O. Zariski \cite{ZO} for curve singularities.
Most of the results were  obtained over the complex numbers.  For example, Bruce and Gaffney \cite{BG} classified the simple parameterized plane curve singularities
over the complex numbers, while Gibson and Hobbs classified the simple parameterizations of space curve singularities, also over the complex numbers \cite{GH1}.
In this article we classify {\em unimodal} (modality 1, see Definition \ref{def.modal}) parametrizations of space curve singularities defined over an algebraically closed field of characteristic $0$. The list of normal forms is given in Tabel \ref{tab:2}.\\

Throughout this paper let $\mathbb{K}$ be an algebraically closed field of characteristic $0$, together with a real valuation $| \, \ | :\K \to \R_{\ge 0}$\,\footnote{This means: $ |a| = 0 \Leftrightarrow a=0$, $ |a||b| = |ab|$, and $|a+b| \le |a|+|b|$. $\K$ is also called a normed field.}, except said otherwise. 
$\K[[x]]$, $x = (x_1,...,x_n)$, denotes the formal power series ring.
For each \mbox{${\eps}\in
  (\R_{>0})^n$}, we define 
   for a formal power series
   \mbox{$H  =  \sum_{{\alpha}\in\N^n}c_{\alpha} {x}^{\alpha} $}, $c_\alpha \in \K$, 
   $\alpha=(\alpha_1,...,\alpha_n)$, 
   $x^\alpha=x_1^{\alpha_1}\cdot...\cdot x_n^{\alpha_n}$,
$$\|H\|_{\eps}:= \sum_{{\alpha} \in \N^n} |c_{\alpha}|\cdot
{\eps}^{{\alpha}}\quad \in \R_{\ge 0}\cup \{\infty\}\,.$$
A formal power series $H$
is called {\em convergent}
iff there exists a real vector \mbox{$\eps\in(\R_{>0})^n$}
such that \mbox{$ \| H \|_{\eps} < \infty $}.
We denote by 
$\K\{x\} \subseteq \K[[x]]$
the {\em convergent  power series ring} over $\K$.  $\K\{x\}/I$, $I$ an ideal,  is a called a local analytic $\K$-algebra (see \cite[Section I.1.2]{GLS25} for properties of analytic $\K$-algebras).
If  $\K$  is complete w.r.t. $| \,\ |$ (every Cauchy sequence with respect to $| \,\ |$ converges in $\K$) then $\K$ is called a {\em complete real valued field}.

\begin{Remark} \label{rm.complete} {\em For a geometric interpretation, we mention that if the valuation $| \,\ |$ on $K$ is non-trivial and complete, then any convergent power series 
 $H= \sum c_\alpha x^\alpha \in \K\{x\}$, $\|F\|_{\eps}< \infty$, defines a continuous  function $H: P^n_\eps \to \K, \ H(a) = \sum c_\alpha a^\alpha$, with $P^n_\eps$ the policylinder
 $$P^n_\eps: =
\{a=(a_1,...,a_n) \in \K^n \mid |a_i| < \eps_i, i =1,...,n\},$$
which is an open neighbourhood (w.r.t the topology defined by the valuation $| \ |$) of $0 \in \K^n$.
The function $H$ is even analytic on $P^n_\eps$, that is, it can be developed in a power series around every point of $P^n_\eps$, as is also the case for $H\in\C\{x\}$ (see \cite[§ 1, Bemerkung]{GR71}). This will be used later.}
\end{Remark}
 
\noindent We call $\K$ an {\em allowed field}, if it is algebraically closed, real valued and of characteristic 0.
\begin{Remark}\label{rm.rvf}{\em
Examples of allowed fields are:
\begin{itemize}
\item Every algebraically closed field $\K$ of characteristic 0 with trivial valuation  ($|0| =0, |a| =1$ if $a\ne 0$), and then $\K\{x\}=\K[[x]]$. 

\item $\K = \C$  with the absolute value as valuation, and  then $\K\{x\}$  is the usual complex analytic power series ring  $\C\{x\}$. Of course, the algebraic closure of $\Q$ is also allowed.
\item The field of {\em Puiseux series}  $\K=\bigcup_{n=1}^\infty F((T^{\frac{1}{n}}))$,  $F$  any algebraically closed field of characteristic $0$, with valuation $|c| = e^{-m/n}$  if $c=\sum_{i=m}^\infty a_iT^{\frac{i}{n}}, a_m\ne0$.
In fact, $\K$ is the algebraic closure of $F((T))$, the quotient field of $F[[T]]$.

\item  Let \mbox{$p$} be a prime number, then the map
$ v_p:\Z\setminus\{0\} \to \R_{>0}\,,\ a \mapsto p^{-m}, \text{ with } m:=\max\{k\in \N  \mid p^k \text{ divides }a\} $
extends to a unique real  {(non-Archimedian)} valuation of $\Q$. The closure of $\Q$ w.r.t. $| \,\ |$ (every Cauchy sequence with respect to $| \,\ |$ converges in the closure) is the field of $p$-adic numbers $\Q_p$. The algebraic closure $\overline \Q_p$ of $\Q_p$ is then an allowed field, as well as $\K = \C_p$, 
the closure of $\Q_p$ w.r.t. $| \,\ |$.
 \end{itemize}}
\end{Remark}

 In this article we want to classify unimodal parametrizations 
\begin{center}
$f : \mathbb{K}\{x,y,z\} \rightarrow t\cdot \mathbb{K}\{t\}, \quad x \mapsto x(t), y \mapsto y(t), z \mapsto z(t)$, 
\end{center}
of space curve singularities,  shortly written as $f=(x(t),y(t),z(t)).$

The following definitions and properties apply more generally to {\em parametrizations} (of reduced, irreducible curve singularities) in $N$-space, that is, $\K$-algebra morphisms 
  $$f : \mathbb{K}\{x_1,...,x_N\} \rightarrow \mathbb{K}\{t\}, \quad x_i \mapsto x_i(t) \in t\cdot \K\{t\},$$
shortly $f= (x_1(t),...,x_N(t))$, such that the {\em$\delta$-invariant} of $f$ is finite,
$$\delta(f) := \dim_\K \K\{t\}/\K\{x_1(t),...,x_N(t)\} < \infty.\footnote{This is equivalent to $f$ being a finite morphism and 
$\K\{x_1,...,x_N\}/\ker(f)  \cong
 \K\{x_1(t),...,x_N(t)\} \subseteq \K\{t\}$ being the analytic local ring of a reduced and irreducible curve singularity (Krull-dimension 1). 
}$$

The finiteness of $\delta$ is equivalent to the existence of a minimal
$c<\infty$ such that the conductor ideal (i.e., the annihilator ideal of ${\mathbb{K}\{t\}}/\K\{x_1(t),...,x_N(t)\}$)
is of the form $t^c\cdot\mathbb{K}\{t\}\subseteq \K\{x_1(t),...,x_N(t)\}.$
The number 
$$c =c(f)$$ 
is called the (multiplicity of) the {\em conductor} of $f$. It is well known that $c(f) \le 2 \delta(f)$ and that equality holds if $f$ is $\mathcal A$-equivalent to a plane curve singularity.
For a  power series $h(t) = \sum_{i\ge m}a_it^i, a_m\ne0,$ we denote the $t$-order of $h$ by 
$\ord_t(h):=m$  and $\ord_t(0):=\infty$, and
define
the {\em order of $f=(x_1(t),...,x_N(t))$}, by 
$$\ord(f) :=  \text{ the minimum of the $t$-orders of $x_1(t),...,x_N(t)$}.$$
Two parametrizations $f= (x_1(t),...,x_N(t))$ and $g= (\hat x_1(t),...,\hat x_N(t))$ are called {\em $\mathcal A$-equivalent} (or right-left equivalent), notation $f\sim g$, if there exist automorphisms 
$$
\begin{array}{lll}
(\psi,\varphi)\in & \mathcal A =\text{Aut}_{\mathbb K}(\mathbb K\{t\})\times\ \text{Aut}_{\mathbb K}(\mathbb 
K\{x_1,...,x_N\}),\\
\psi : &\mathbb{K}\{t\}\rightarrow \mathbb{K}\{t\},\\
\varphi : &\mathbb{K}\{x_1,...,x_N\}\rightarrow 
\mathbb{K}\{x_1,...,x_N\},
\end{array}
$$
such that the following diagram commutes:
\vskip 0,65 true cm
\begin{center}
		\begin{picture}(300,70)
		\put(30,54){\vector(0,-1){33}} \put(248,54){\vector(0,-1){33}}
		\multiput(62,12)(0,50){2}{\vector(1,0){76}}
		\put(0,10){$\K\{x_1,...,x_N\}$} \put(140,10){$\K\{\hat x_1(t),...,\hat x_N/t)\}\subseteq \K\{t\}$}
		\put(0,60){$\K\{x_1,...,x_N\}$} \put(140,60){$\K\{x_1(t),...,x_N(t)\}\subseteq \K\{t\}$}
		\put(1,35){$\varphi$} \put(254,35){$\psi$}
		\put(80, 15){$g$} \put(80, 65){$f$}
		\end{picture}
	\end{center}
that is, the (right-left) group $\mathcal A$ acts by 
$\mathcal A \cdot f = g = \psi \circ f \circ \varphi^{-1}.$

\begin{Remark}\label{rm.aut}{\em
To show that $g$ is in the $\mathcal A$-orbit of $f$, we have to show the existence of $\varphi$ and $\psi$ satisfying 
$\psi \circ f = g\circ \varphi$. That is, if
 $\varphi(x_i)=H_i(x_1,...,x_N)$, then
$$(x_1(\psi(t)),...,x_N(\psi(t)))=
(H_1(\hat x_1(t),...,\hat x_N(t)),
...,H_N(\hat x_1(t),...,\hat x_N(t))).$$
An automorphism $\psi \in \text{Aut}_{\K}(\K\{t\})$ is of the form 
$\psi(t) = a_1t+ a_2t+a_3t +....$, with $a_1\ne 0$ and $a_i\in \K$ arbitrary for $i>1$. Setting $\alpha = a_1^{-1}$, we can replace $\psi(t)$ by $\psi'(t) = \psi(\alpha t) = t + a'_2t+a'_3t +....$  
and get as equivalent equation 
$$(x_1(\psi'(t)),...,x_N(\psi'(t)))=
(H_1(\hat x_1(\alpha t),...,\hat x_N(\alpha t)),
...,H_N(\hat x_1(\alpha t),...,\hat x_N(\alpha t))).$$
We will use this form later, as it is computationally much faster.
}
\end{Remark}

\begin{Definition}\label{def.gamma}
Let $f=(x_1(t),...,x_N(t))$ be a parametrization of a curve singularity.
\begin{enumerate}
\item  We define the {\em semigroup} of $f$ as
\begin{center}
$\Gamma={\Gamma}(f)= \{{{\ord}_t(h)} \mid h\in {\mathbb{K}\{x_1(t),...,x_N(t)\}\setminus \{0\}}\}$.
\end{center}
\item If $t^c\cdot\mathbb{K}\{t\}$ is the conductor
ideal then $c-1\not \in \Gamma$ and $l\in \Gamma$ if
$l\geq c$. The integer $c=c(\Gamma)$ is called the {\em conductor} of \ $\Gamma$ or of $f$. 
\item The cardinality of the set $\mathbb Z_{\geq 0} \setminus \Gamma$ is called the {\em $\delta$-invariant} and denoted by
$\delta=\delta({\Gamma}(f))$. It coincides with $\delta(f)$.
\item The semigroup $\Gamma$ has a unique {\em minimal system of generators} $\{\beta_1,\ldots,\beta_k\}$ and we write $\Gamma=\langle  \beta_1,\ldots,\beta_k \rangle$. We will always assume that the
minimal generators of a semigroup are given in an increasing way, 
$$\beta_1 < \beta_2 < ... < \beta_n.$$
We call $\beta_1$ the {\em order} of $\Gamma$.
\end{enumerate}
\end{Definition}

In order to define the modality of a parametrization $f$, we have to consider deformations over base spaces, which have {\em closed points} in the neighbourhood of a given point.
This is the case for the spectrum of $\K[s]$, however not for e.g.  $\K[[s]]$, $s=(s_1,...,s_n)$. Therefore we will use affine algebras $A \cong \K[s]/I$ as base algebras for deformations of a parametrization. The base space is then $\Max A$ (the set of maximal ideals of $\Spec A$), identified with the closed points of the algebraic subvariety $V(I) \subseteq \K^n$, endowed with the Zariski topology.  
\medskip

In the complex analytic situation, the base algebra of a deformation is usually of the form $\C\{s\}/I$. The base space is then represented by a small neighborhood  of the origin in the analytic set $V(I)\subseteq \C^n$, with the (Euclidean) $\C$-topology. 
The same holds, if the valuation $| \ |$ on $\K$ is complete and non-trivial. Then a convergent power series can be identified with an analytic funtcion on an open neighbourhood of $0\in \K^n$ (e.g. on $P_\eps$ for some $\eps >0$, see Remark \ref{rm.complete}). Hence, if we use $\K\{s\}/I$ as base algebra, the base space is a small neighborhood of the origin of the analytic set $V(I)\subseteq \K^n$, 
with the $\K$-topology, i.e., the topology on $\K$ is given by the valuation (as for $\K=\C$).\\

In order to be able to
treat deformations of parametrizations of e.g. formal power series over Zariski open sets and analytic power series over $\K$-open sets at the same time, we make the following definition. 
\begin{Definition} \label{def.base} Let $\K$ be any allowed field, $s=(s_1,...,s_n)$.
\begin{enumerate} 
\item We call a $\K$-algebra $A$ {\em geometric}, if, for some proper  ideal $I$, $A$ is\\
$(i)$ isomorphic to $\K[s]/I$, or\\
$(ii)$ isomorphic to $\K\{s\}/I$ and the valuation on $\K$ is non-trivial and complete.
\item If  $A \cong \K[s]/I$ in case (i), we
denote by $B_A$ the set of maximal ideals of $\Spec A$ containing $I$, which we identify with the closed points of $V(I) \subseteq \K^n$. 

\item If  $A \cong \K\{s\}/I$ in case (ii), there exists an open neighbourhood $U$ of $0\in \K^n$ such that any finitely many generators of $I$ can be considered  as analytic functions on $U$ (by Remark \ref{rm.complete}) and the analytic set $V(I) \subseteq U $ is defined. In this case we denote by $B_A$ the set $V(I)$ (which depends of course on $U$). 

\item If $p \in B_A$ and 
$\fm_p$  the maximal ideal corresponding to $p$, then  $\K = A/\fm_p$ and we have a natural map $A[[t]] \to \K[[t]]$,  $H(s,t) \mapsto H(p,t)$. We set
$$A\{t\} := \{H(s,t) \in A[[t]] \mid H(p,t) \in \K\{t\} \text{ for } p \in B_A\}.$$ 
\end{enumerate}
 We call $B_A$ the {\em geometric  space} of the geometric algebra $A$.
$B_A$ will be the base space for our families of parametrizations. It carries the Zariski topology in case 2. and the $\K$-topology in case 3. and consists of closed points of $\K^n$.\footnote{In the following a point means always a closed point.} While case 3. is only used if the valuation  on $\K$ is complete, case 2. will be used  for any allowed field.
\end{Definition}

\begin{Definition}[Family of parametrizations]\label{def.unf} 
Let $A$ be a geometric algebra and $B_A$ its geometric space. A {\em family of parameterizations} over $B_A$ (or over $A$) is an element $F=(X_1,...,X_N) \in A\{t\}^N$ such that 
for each $p\in B_A$ 
$$F(p,t) = (X_1(p,t),..., X_N(p,t)) \in t\cdot \K\{t\}^N$$
  is a parametrization of a  curve singularity. Instead of $F$ we write more intuitively $F(s,t)$  or $(X_1(s,t),..., X_N(s,t))$. 

For fixed $p$ and any open neighbourhood $U$ of $p$ in $B_A$ the restriction of the family $F$ to $U$ is then called an {\em unfolding} (or {\em deformation with trivial section)} of the parametrization $F(p,t)$. 
\end{Definition} 

Recall that a parametrization $f=(x_1(t),...,x_N(t))\in t\cdot {\K\{t\}}^N$ is {\em $k$-determined} (w.r.t. $\mathcal A$), if, for every parametrization $g \in t\cdot {\K\{t\}}^N$, $f \equiv g$ mod $t^{k+1} {\K\{t\}}^N$ implies $f\sim g$.
By definition,  the delta invariant  $\delta(f)$  and hence $c(f)$ for a parametrization  $f$ is finite. By \cite [Proposition 11]{GP}  
$f$ is $(2c-1)$ and hence $(4\delta -1)$-determined (the proof in loc. cit.  is given for formal power series, but it works as well in the convergent case).

 Denote by $f^{(k)}$ the $k$-jet of $f$, that is the image of $f$ in ${\K\{t\}}^N / t^{k+1} {\K\{t\}}^N$, identified with power series expansions of $x_1(t),...,x_N(t)$ up to (including) order $k$. Then the $k$-jet $\mathcal A^{(k)}$ of $\mathcal A$ is an algebraic group, which 
 acts on the affine space ${\K\{t\}}^N / t^{k+1} {\K\{t\}}^N$ by letting $\mathcal A^{(k)} \cdot f^{(k)}$ be the $k$-jet of $\mathcal A \cdot f$.  Note that  $f^{(k)} \sim f$ if $k\ge 2c-1$.\\

In some cases the determinacy bound can be improved.
 \begin{Lemma} \label{lem.determined}
 Let $f=(x_1(t),...,x_N(t))\in t\cdot {\K\{t\}}^N$ be a parametrization and $\Gamma$ be the semigroup with conductor $c$ and minimal set of generators $\beta_1,\ldots,\beta_n$. If $\beta_n\leq c$ then $f$ is 
 $c$-determined and if $\beta_n < c$ then $f$ is 
 $(c-1)$-determined.
 \end{Lemma}
 \begin{proof}
Let $(\bar x_1(t),...,\bar x_N(t)):=f^{(c)}$  be the 
$c$-jet of $f$. Then, for each $i=1,...,N$, there exists $g_{i}(t) \in \K\{t\}$ such that 
$$x_i(t)= \bar x_i(t) +g_i(t), \text{ with }\ord(g_i) \geq c+1.$$
Since $g_i$ is contained in the conductor ideal of the parametrization $f$, we have 
$$ \K\{\bar x_1(t),...,\bar x_N(t)\} \subseteq 
\K\{x_1(t),...,x_N(t)\}.$$
In particular $\bar \Gamma \subseteq \Gamma$, where $\bar\Gamma$ is the semigroup of the parametrization $(\bar x_1(t),...,\bar x_N(t))$.

 Now let $h_j(t) \in \K\{x_1(t),...,x_N(t)\}$ be such that 
 $\beta_j = \ord(h_j)$, $j=1,...,n$. Then, for each $j$,  there exist 
 $a_j^i \in \K\{t\}$, $i=1,...,N$,  such that
$h_j(t) = \sum _{i=1}^Na_j^i x_i(t)$ for $j=1,...,n$. 
Taking $c$-jets, we get
$$\bar h_j(t) = \sum _{i=1}^N \bar a_j^i \bar x_i(t), \ j=1,...,n.$$
Since $\beta_n \leq c$ we have $\beta_j = \ord(\bar h_j(t))$ and hence $\beta_j \in \bar \Gamma$. This shows $\bar \Gamma = \Gamma$ and that $c$ is the conductor of $(\bar x_1(t),...,\bar x_N(t))$. Therefore $ \K\{x_1(t),..., x_N(t)\} \subseteq 
\K\{\bar x_1(t),...,\bar x_N(t)\}$ and with the other inclusion
$ \K\{\bar x_1(t),...,\bar x_N(t)\} = 
\K\{x_1(t),...,x_N(t)\}$. 

If $\beta_n < c$ we can argue in the same way to get that $f$ is 
 $(c-1)$-determined.
This proves the lemma.
 \end{proof}

Let us start with an informal definition of the modality of a parametrization. 

\begin{Remark}\label{rm.modal}{\em 
Let $f=(x_1(t),...,x_N(t))\in t\cdot {\K\{t\}}^N$ define a parametrized curve singularity. Then the {\em modality of $f$ is $\ge n$} (or $f$ is {\em at least $n$-modal}), if there exists an unfolding  $F(s,t) \in t\cdot A\{t\}^N$ of $f$ and an open subset $U \subseteq B_A$  containing $0$, such that the set
$\{F(p,t)\in \mathbb{K}\{t\}^N \mid p\in U \}$ 
of parametrizations contains  
$n$-dimensional families of $\mathcal A$-equivalence classes. 
 $f$ is called {\em simple} (resp. {\em unimodal}, resp. {\em bimodal}) if the modality of $f$ is $0$ (resp. $1$, resp. $2$).}
\end{Remark}

The informal definition of the modality is in principal sufficient for computations. However, we need also the semicontinuity of the modality. This can be proven using the  (properly modified) formal definition of the modality from \cite{GH}, which we will now outline.

 For the definition of the modality of a parametrization $f$ we use the modality w.r.t. the action of an algebraic group $G$ on an algebraic variety $X$, as defined in \cite[Definition 2.1]{GH}.

\begin{Definition}[Modality]\label{def.modal}
\begin{enumerate}
\item Let $\{X_{i}, i=1,\ldots,s\}$ be a Rosenlicht stratification\footnote{This is a decomposition of $X$ as finite disjoint union of $G$-invariant locally closed algebraic subvarieties $X_{i}$, such that $X_{i}/G$ is a geometric quotient with quotient morphism $p_{i}\colon X_{i}\to X_{i}/G$.} of the algebraic variety $X$ under the action of an algebraic group $G$ with quotient morphisms $p_i: X_i\to X_i/G$, and let $U$ be an open neighbourhood of $x\in X$. We define
$$G\text{-}\mathrm{mod}(U):=\max_{1\leq i\leq s}\{\dim \big (p_i(U\cap X_i)\big)\},$$
and call 
$$G\text{-}\mathrm{mod}(x):=\min \{G\text{-}\mathrm{mod}(U)\ |\ U \text{ a neighbourhood of x}\}$$
the {\em $G$-modality} of $x$ (in $X$).

\item The {\em modality of the parametrization $f$}, denoted by {\em mod}$(f) =$ \mbox{$\mathcal A$-{\em mod}$(f)$},  is defined as the $\mathcal A^{(k)}$-modality of $f^{(k)}$ under the action of $\mathcal A^{(k)}$ on ${\K\{t\}}^N / t^{k+1} {\K\{t\}}^N$ for $k\ge 4\delta-1$ (which is independent of the choice of $k\ge 4\delta-1$, see \cite[Proposition 2.6]{GH}).
\end{enumerate}
\end{Definition}
 
 \begin{Proposition}\label{prop.sem}
 The modality of a parametrized curve singularity is upper semicontinuous for unfoldings. That is, for any unfolding  
$F\in t\cdot {A\{t\}}^N$ of $f=F(0,t)$ there exists a neighbourhood $U \subseteq B_A$ of $0$, such that 
 {\em mod}$(f)$ $\ge$  {\em mod}$(F(p,t))$ for all $p\in U$.\\
 In particular, if  {\em mod}$(F(p,t))$ $\ge 2$ for all  $p\in U \setminus Z$, $Z$ a nowhere dense subset of $U$
 then $f$ is not simple and not unimodal.
 \end{Proposition}
 
 \begin{proof} The proof is basically the same as the proof of Proposition 2.7 in \cite{GH}. We use the semicontinuity of 
 $\delta$ for parametrizations (\cite[Theorem 18]{GP}), which implies the simultaneous bound   $4\delta-1$, $\delta = \delta(f)$,  for the determinacy  of $F(p,t)$, $p\in U$. It follows from Definition 
 \ref{def.modal},  that  \mbox{$\mathcal A^{(k)}$-{mod}$(F(p,t)^{(k)})$} is semicontinuous. For $k\ge 4\delta-1$ we have $F(p,t)^{(k)} \sim F(p,t)$, which implies the result.
 \end{proof}
 
\begin{Remark}\label{rm.smod}{\em
 The modality of a parametrization $f:\K\{x_1,...,x_N\} \to \K\{t\}$  is not an invariant of the abstract curve singularity  
 $\K\{x_1(t),...,x_N(t)\} \cong
 \K\{x_1,...,x_N\}/\ker (f)$, but depends on the embedding in $N$-space. Consider for example 
 $f=(t^5,t^7)$ and $\tilde f = (t^5,t^7,0)$.  Then $\K\{x_1(t),x_2(t)\}\cong \K\{x_1(t),x_2(t),x_3(t)\}$, but by \cite[Theorem 2]{MP} $f$ has modality $1$ while the modality of $\tilde f$ is $\ge 2$ by Proposition \ref{1.3.1}.}
 \end{Remark}
 \begin{Definition}[Modular family]\label{def.modular} 
Let  $F(s,t)$  be a family of parametrizations over $B_A$.
$F$ is called {\em ($\mathcal A$-)modular} if there exists a nowhere dense subset $Z \subset B_A$ such that $F(p,t)$ and $F(q,t)$ are not $\mathcal A$-equivalent 
for each $p\ne q$ in $B_A \setminus Z$. 
\end{Definition}
Informally speaking, a family of parameterizations is modular if 
the parametrizations at two distinct generic points  of every irreducible component of $B_A$ are not $\mathcal A$-equivalent. The following proposition shows that the modality is bounded from below by the dimension of some modular family. This is quite useful and will be used in our classification. The proof of the proposition follows easily 
from the semicontinuity of the modality (Proposition \ref{prop.sem}).
(see also \cite[Proposition 2.10]{BGMP} for
parametrizations of plane curve singularities with a slightly weaker notion of 'modular family').
\begin{Proposition}\label{prop.modular}
Let $f$ be a parameterization of a space curve singularity. Moreover, let $F(s,t)$ be a family of parameterizations over $B_A$  and $p \in B_A$ such that $f \sim F(p,t)$. 
If $F$ is modular, then
$\mathcal A\text{\em -mod}(f) \geq \dim (B_A,p).$\footnote{$\dim (B_A,p)$ is the dimension of the local ring of $B_A$ at $p$, which is equal to $\dim B_A$ if $B_A$ is irreducible.}
\end{Proposition}

We now recall the classification of simple parametrized space curve singularities due to Gibson and Hobbs \cite{GH1} for complex analytic parametrizations.

\begin{Theorem}\label{ch0}
Let $\mathbb K$ be any allowed field and
 $f =(x(t),y(t),z(t)) \in t \cdot {\mathbb{K}\{t\}}^3$ a parametrized space curve singularity.  Then $f$ is simple iff it is $\mathcal A$-equivalent to a parametrized space curve singularity in Table 1.
\end{Theorem}

\begin{table}[htbp]
\centering
\renewcommand{\arraystretch}{1.15}
\setlength{\tabcolsep}{4pt}
\caption{}
\label{tab:GH}
\begin{tabular}{|c|c|l|l|}
\hline
$\Gamma(f)$ & $c(f)$ & Normal Form & Condition\\
\hline
$\langle 1 \rangle$ & $0$ & $(t, 0, 0)$ & \\
\hline
$\langle 2, 2k+1 \rangle$ & $2k$ & $(t^2, t^{2k+1}, 0)$ & $k \geqslant 1$ \\
\hline
$\langle 3, 3k+1 \rangle$ & $6k$ &
\stack{$(t^3, t^{3k+1}, 0)$ \\ $(t^3, t^{3k+1}+t^{3p+2}, 0)$} &
\stack{$k \geqslant 1$ \\ $k \leqslant p < 2k-1$} \\
\hline
$\langle 3, 3k+1, 3n+2 \rangle$ & $3n$ &
\stack{$(t^3, t^{3k+1}, t^{3n+2})$ \\ $(t^3, t^{3k+1}+t^{3p+2}, t^{3n+2})$} &
\stack{$k \leqslant n < 2k$ \\ $k \leqslant p < n < 2k$} \\
\hline
$\langle 3, 3k+2 \rangle$ & $6k+2$ & $(t^3, t^{3k+2}, 0)$ & $k \geqslant 1$ \\
\hline
$\langle 3, 3p+2 \rangle$ & $6p+2$ & $(t^3, t^{3k+1}+t^{3p+2}, 0)$ & $p < k < 2p$ \\
\hline
$\langle 3, 3n+2, 3k+1 \rangle$ & $3k-1$ &
\stack{$(t^3, t^{3k+1}, t^{3n+2})$ \\ $(t^3, t^{3n+2}+t^{3p+1}, t^{3k+1})$} &
\stack{$n < k \leqslant 2n$ \\ $n < p < k \leqslant 2n$} \\
\hline
$\langle 4, 5 \rangle$ & $12$ &
\stack{$(t^4, t^5, 0)$ \\ $(t^4, t^5+t^7, 0)$} & \\
\hline
$\langle 4, 5, 6 \rangle$ & $8$ & $(t^4, t^5, t^6)$ & \\
\hline
$\langle 4, 5, 7 \rangle$ & $7$ & $(t^4, t^5, t^7)$ & \\
\hline
$\langle 4, 5, 11 \rangle$ & $8$ &
\stack{$(t^4, t^5, t^{11})$ \\ $(t^4, t^5+t^7, t^{11})$} & \\
\hline
$\langle 4, 6, 2k-1, 2k+1 \rangle$ & $2k-2$ & $(t^4, t^6+t^{2k-7}, t^{2k+1})$ & $k \geqslant 7$ \\
\hline
$\langle 4, 6, 2k+1 \rangle$ & $2k+4$ &
\stack{$(t^4, t^6, t^{2k+1})$ \\ $(t^4, t^6+t^{2k-1}, t^{2k+1})$} &
\stack{$k \geqslant 3$ \\ $k \geqslant 4$} \\
\hline
$\langle 4, 6, 2k+1, 2k+3 \rangle$ & $2k$ & $(t^4, t^6+t^{2k-3}, t^{2k+1})$ & $k \geqslant 5$ \\
\hline
$\langle 4, 6, 2k+7 \rangle$ & $2k+10$ & $(t^4, t^6+t^{2k+1}, 0)$ & $k \geqslant 3$ \\
\hline
$\langle 4, 7 \rangle$ & $18$ &
\stack{$(t^4, t^7, 0)$ \\ $(t^4, t^7+t^9, 0)$ \\ $(t^4, t^7+t^{10}, 0)$} & \\
\hline
$\langle 4, 7, 9 \rangle$ & $11$ &
\stack{$(t^4, t^7, t^9+t^{10})$ \\ $(t^4, t^7, t^9)$} & \\
\hline
$\langle 4, 7, 10 \rangle$ & $14$ &
\stack{$(t^4, t^7, t^{10})$ \\ $(t^4, t^7+t^9, t^{10})$} & \\
\hline
$\langle 4, 7, 13 \rangle$ & $11$ &
\stack{$(t^4, t^7, t^{13})$ \\ $(t^4, t^7+t^9, t^{13})$} & \\
\hline
$\langle 4, 7, 17 \rangle$ & $14$ &
\stack{$(t^4, t^7, t^{17})$ \\ $(t^4, t^7+t^9, t^{17})$ \\ $(t^4, t^7+t^{10}, t^{17})$} & \\
\hline
\end{tabular}
\end{table}


\begin{Remark}{\em
(1) The theorem was proven by Gibson and Hobbs (see \cite{GH1}) for $\K=\mathbb{C}$, the field of complex numbers.
In \cite{JP}, this classification was also carried out for algebraically closed fields of characteristic $p>0$. If  $p>13$, the results of \cite{JP} and \cite{GH1} agree, and the methods used in \cite{JP} also work for characteristic zero for every allowed field. This proves the theorem.

(2) Compared to  \cite[Table 1]{GH1} we changed in line   $\langle 3,3k+1\rangle$ the condition 
$k \le p < 2k$ to $k \le p < 2k-1$ and in line $\langle 3,3p+2\rangle$ the condition $p < k \le 2p$ to $p < k < 2p$.
In both cases the term $t^{c-1}$ can be deleted due to Proposition \ref{prop.Za}.
}
\end{Remark}

\begin{Corollary} \label{cor.cho}
Every parametrization of a space curve $f=(x(t),y(t), z(t))$ with semigroup $\Gamma (f)\leq \langle 4,7\rangle$ (see Definition \ref{1.2.1}) is simple.
\end{Corollary}
\begin{proof}
By inspection of Table \ref{tab:GH}, every parametrization $f$ of order $3$ or with $\Gamma(f)$ containing $4,5$ and $4,6$ is simple. It remains to consider semigroups of order $4$ and containing $4,7$. These are, besides those of Table \ref{tab:GH}, the
semigroups
$\langle 4,7,9,10\rangle$ and $\langle 4,7,10,13\rangle$.
But these cannot be realized by a space curve, which can easily be checked.
\end{proof}

An important  basic step for the classification is a theorem of Zariski (\cite{ZO}, Chapter III, Proposition 1.2) on "short parametrizations", which can be generalized to space curves. We omit the proof, as it is similar to Zariski's.
\begin{Theorem}\label{thm.Za}
Let $\mathbb K$ be an algebraically closed field of characteristic $0$.
Given a parametrization
\begin{center}
 $(t^l+{\Sigma}_{i>l} \ a_it^i, \ t^m+{\Sigma}_{i>m}\ b_it^i, \ t^n+ {\Sigma}_{i>n}\ c_it^i)$
 \end{center}
 $2<l<m<n$, $l\nmid m$, and $n \notin \langle l,m\rangle$ if $n<\infty$,
 with semigroup $\Gamma$ and  conductor $c$. For $k\in \Gamma$  there exists an $\mathcal A$-equivalent parametrization
 \begin{center}
 $(t^l+{\Sigma}_{i>l}\ \hat{a}_it^i, \ t^m+{\Sigma}_{i>m}\ \hat{b}_it^i, \ t^n + {\Sigma}_{i>n}\ \hat{c}_it^i)$
 \end{center}
 with $a_i={\hat{a}}_i,b_i={\hat{b}}_i, c_i={\hat{c}}_i $ if $i<k$, ${\hat{a}}_k={\hat{b}}_k={\hat{c}}_k=0$, and ${\hat{a}}_s={\hat{b}}_s={\hat{c}}_s=0$,  for all $s\geq c$.
\end{Theorem}

The following proposition is an  improvement of Zariski's short parametrization. Together wir Remark \ref{rm.Za}, it is important for the classification.

\begin{Proposition}\label{prop.Za}
Let
$(x(t),y(t),z(t)) \in t\cdot \K\{t\}^3$ be a parametrization with semigroup $\Gamma$ and conductor $c$.
\begin{enumerate}  
\item Then there exists an $\mathcal A$-equivalent parametrization of the form
 $$f=(t^l, \ t^m+{\Sigma}_{i>m,i \notin \Gamma}\ {b}_it^i, \ t^n + {\Sigma}_{i>n,i \notin \Gamma} \ {c}_it^i),$$
$l<m<n$ ($n=\infty$ included), $l\nmid m$ and $n \notin \langle l,m \rangle $ if $n<\infty$.
\item If $c-1-m+l\in \Gamma$ and $c-1-n+l\in \Gamma$ then the parametrization $f$ is $\mathcal A$-equivalent to a parametrization of the form
$$(t^l, \ t^m+{\Sigma}_{c-1>i>m,i \notin \Gamma}\ {\hat b}_it^i, \ t^n + {\Sigma}_{c-1>i>n,i \notin \Gamma} \ {c}_it^i)$$
for certain $\hat b_i \in \K$ (with the same $c_i$ as in $f$).\\
(a) If only $c-1-n+l\in \Gamma$ then the parametrization $f$ is 
$\mathcal A$-equivalent to $(t^l, \ t^m+{\Sigma}_{i>m,i \notin \Gamma}\ {\hat b}_it^i, \ t^n + {\Sigma}_{c-1>i>n,i \notin \Gamma} \ {c}_it^i)$\\
(b) If  only $c-1-m+l\in \Gamma$ then the parametrization $f$ is $\mathcal A$-equivalent to $(t^l, \ t^m+{\Sigma}_{c-1>i>m,i \notin \Gamma}\ {b}_it^i, \ t^n + {\Sigma}_{i>n,i \notin \Gamma} \ {c}_it^i)$. 
\end{enumerate}
\end{Proposition}

In addition to Proposition \ref{prop.Za} we will use also the following remark to simplify a parametrization. It says that we can make the coefficient of the 2nd non-vanishing term of either $y(t)$ or $z(t)$ equal to 1 (but not of $y$ and $z$ simultaneously).

\begin{Remark}\label{rm.Za} {\em
Let 
$f=(x(t),y(t),z(t))=(t^l, \ t^m+\alpha_1 t^{m_1} + \alpha_2 t^{m_2}  +...,
 \ t^n +\beta_1 t^{n_1} + \beta_2 t^{n_2}  +...)$
be a parametrization with $m<m_1<m_2<...$, $n<n_1<n_2<...$ and $\alpha_i, \beta_i \in \K\setminus \{0\}$.\\
We use the map $$t\mapsto \beta^{-1}t, \ \beta ^{m_1-m}
= \alpha_1$$ 
and multiply the components with  $ \beta^{l}, \beta^{m}, \beta^n$ respectively,
to obtain
$$f\sim (t^l, \ t^m+t^{m_1} + \tilde\alpha_2 t^{m_2} + ...,
 \ t^n +\tilde\beta_1 t^{n_1} + \tilde\beta_2 t^{n_2}  +...).
$$ 
Moreover, $\tilde \beta_i = \beta_i \alpha_1^{(m_1-m)(n_i-n)}$ covers all elements of $K$ if $\alpha_1$ does.
The same argument works for the third component $z(t)$.}
\end{Remark}

\begin{proof}[Proof  of the Proposition]
The first part is an immediate consequence of Theorem \ref{thm.Za}.\\
For 2. assume that $c-1-m+l\in \Gamma$ and 
 $c-1-n+l\in \Gamma$.  To construct the normal form, we have to start with the 3rd component.\\
We use first that $c-1-n+l\in \Gamma$ and define an automorphism of $\K\{t\}$ by
$$\varphi(t):=t-\frac{1}{n}c_{c-1}t^{c-n}.$$
Then $$\varphi(t^m+{\Sigma}_{i>m,i \notin \Gamma}\ {b}_it^i)=t^m+{\Sigma}_{i>m}\ {\bar b}_it^i$$
for suitable $\bar b_i$,
$$\varphi(t^n + {\Sigma}_{i>n,i \notin \Gamma} \ {c}_it^i)=t^n + {\Sigma}_{c-1>i>n,i \notin \Gamma} \ {c}_it^i+\text{terms of order}\geq c, \text{ and} $$
$$\varphi(t^l)=t^l-\frac{lc_{c-1}}{n}t^{c-1-n+l}+\text{higher order terms}.$$
Since $c-1-n+l\in \Gamma$, using Theorem \ref{thm.Za}, we can kill the terms with order in $\Gamma$ and especially the term 
$\frac{lc_{c-1}}{n}t^{c-1-n+l}$. We obtain that
$$\varphi((t^l, \ t^m+{\Sigma}_{i>m,i \notin \Gamma}\ {b}_it^i, \ t^n + {\Sigma}_{i>n,i \notin \Gamma} \ {c}_it^i))$$
is $\mathcal A$-equivalent to
$$(t^l+\text{terms of order} > c-1-n+l,t^m+{\Sigma}_{i>m}\ {\bar b}_it^i,t^n + {\Sigma}_{c-1>i>n,i \notin \Gamma} \ {c}_it^i).$$
Now choose the automorphism $\psi$ with the property $\psi(t^l+\text{terms of order} > c-1-n+l)=t^l$.
Then 
$$\psi((t^l+\text{terms of order} > c-1-m+l,t^m+{\Sigma}_{i>m}\ {\bar b}_it^i,t^n + {\Sigma}_{c-1>i>n,i \notin \Gamma} \ {c}_it^i))$$
$$=(t^l,t^m+{\Sigma}_{i>m}\ {\tilde b}_it^i,t^n + {\Sigma}_{c-1>i>n,i \notin \Gamma} \ {c}_it^i+\text{terms of order} \geq c)$$
for suitable $\tilde b_i$.
Again using Theorem \ref{thm.Za} we can kill terms of the second component and obtain that 
$(t^l, \ t^m+{\Sigma}_{i>m,i \notin \Gamma} \ {b}_it^i, \ t^n + {\Sigma}_{i>n,i \notin \Gamma} \ {c}_it^i),$
is $\mathcal A$-equivalent to
$$(t^l,t^m+{\Sigma}_{i>m,i \notin \Gamma}\ {\hat b}_it^i,t^n + {\Sigma}_{c-1>i>n,i \notin \Gamma} \ {c}_it^i)$$
for suitable $\hat b_i$. This proves already 2.(a).\\
Now assume that also $c-1-m+l\in \Gamma$. We repeat the same procedure for the second component of the parametrization. Let
$$\varphi(t):=t-\frac{1}{m}\hat b_{c-1}t^{c-m}.$$
Then 
$$\varphi(t^m+{\Sigma}_{i>m,i \notin \Gamma}\ {\hat b}_it^i)=t^m+{\Sigma}_{c-1>i>m,i \notin \Gamma}\ {\hat b}_it^i+\text{terms of order}\geq c$$
$$\varphi(t^n + {\Sigma}_{c-1>i>n,i \notin \Gamma} \ {c}_it^i)=t^n + {\Sigma}_{c-1>i>n,i \notin \Gamma} \ {c}_it^i+\text{terms of order}\geq c$$
$$\varphi(t^l)=t^l-\frac{l\hat b_{c-1}}{m}t^{c-1-m+l}+\text{higher order terms}.$$
Since $c-1-m+l\in \Gamma$, using Theorem \ref{thm.Za}, we can kill the terms of order $\geq c$ and the term $\frac{l\hat b_{c-1}}{m}t^{c-1-m+l}$ and obtain that
$$\varphi((t^l, \ t^m+{\Sigma}_{i>m,i \notin \Gamma}\ {\hat b}_it^i, \ t^n + {\Sigma}_{c-1>i>n,i \notin \Gamma} \ {c}_it^i))$$
is $\mathcal A$-equivalent to
$$(t^l+\text{terms of order} > c-1-m+l,t^m+{\Sigma}_{c-1>i>m,i \notin \Gamma}\ {\hat b}_it^i,t^n + {\Sigma}_{c-1>i>n,i \notin \Gamma} \ {c}_it^i).$$
As before we define $\psi$ with the property $\psi(t^l+\text{terms of order} > c-1-m+l)=t^l$ and obtain
\begin{center}
$\psi((t^l+\text{terms of order} > c-1-m+l,t^m+{\Sigma}_{c-1>i>m,i \notin \Gamma}\ {\hat b}_it^i,t^n + {\Sigma}_{c-1>i>n,i \notin \Gamma} \ {c}_it^i))=$
$(t^l,t^m+{\Sigma}_{c-1>i>m,i \notin \Gamma}\ {\hat b}_it^i+\text{terms of order} \geq c,t^n + {\Sigma}_{c-1>i>n,i \notin \Gamma} \ {c}_it^i+\text{terms of order} \geq c).$
\end{center}
Again using Theorem \ref{thm.Za} we obtain that 
$(t^l, \ t^m+{\Sigma}_{i>m,i \notin \Gamma}\ {\hat b}_it^i, \ t^n + {\Sigma}_{c-1>i>n,i \notin \Gamma} \ {c}_it^i),$
is $\mathcal A$-equivalent to
$$(t^l,t^m+{\Sigma}_{c-1>i>m,i \notin \Gamma}\ {\hat b}_it^i,t^n + {\Sigma}_{c-1>i>n,i \notin \Gamma} \ {c}_it^i).$$
 If only  $c-1-m+l\in \Gamma$, the same argument (with $b_i$ instead of $\hat b_i$) proves the claim of 2.(b).
\end{proof}

\begin{Definition}
A parametrization $f=(x(t),y(t),z(t))$ with the properties of  Proposition \ref{prop.Za} is called a {\em short parametrization}.
\end{Definition}

The main result of this paper is the following theorem:
\begin{Theorem}\label{101}

  Let $f = (x(t), y(t), z(t))$ be a parametrization of a unimodal space
  curve then it is $\mathcal A$-equivalent to a parametrization in Table 2.\footnote{No condition on $a$ means $a\in \K$; modulus $64a^2-63a-56b$ means that the equivalence classes in the $(a,b)$-plane are the level sets $64a^2-63a-56b =$ constant.}
  \end{Theorem}

\begin{table}[htbp]
\centering
\renewcommand{\arraystretch}{1.2}
\setlength{\tabcolsep}{5pt}
\caption{}
\label{tab:2}
\begin{tabular}{|c|c|l|l|}
\hline
$\Gamma(f)$ & $c(f)$ &  \multicolumn{1}{c|}{Normal Form} & Condition\\
\hline
$\langle 4, 9, 10 \rangle$ & $16$ &
\stack{$(t^4, t^9+t^{11}, t^{10}+at^{11})$ \\ 
$(t^4, t^9, t^{10}+at^{11})$} &
\stack{\strut \\ $a \in \{0,1\}$} \\
\hline
$\langle 4, 9, 11 \rangle$ & $15$ &
\stack{$(t^4, t^9+t^{10}, t^{11}+at^{14})$ \\ $(t^4, t^9, t^{11}+at^{14})$} &
\stack{\strut \\ $a \in \{0,1\}$} \\
\hline
$\langle 5, 6, 7 \rangle$ & $10$ &
\stack{$(t^5, t^6, t^7+t^8+at^9)$ \\ $(t^5, t^6, t^7+at^9)$} &
\stack{\strut \\$a \in \{0,1\}$} \\
\hline
$\langle 5, 6, 8 \rangle$ & $10$ &
\stack{$(t^5, t^6+at^7, t^8)$} &
\stack{$a \in \{0,1\}$} \\
\hline
$\langle 5, 6, 9 \rangle$ & $14$ &
\stack{$(t^5, t^6+t^7+at^8, t^9)$ \\ $(t^5, t^6+at^8, t^9)$} &
\stack{\strut \\ $a \in \{0,1\}$} \\
\hline
$\langle 5, 7, 8 \rangle$ & $12$ &
\stack{$(t^5, t^7+t^9+at^{11}, t^8)$ \\ $(t^5, t^7+at^{11}, t^8)$ \\ $(t^5, t^7+at^9+bt^{11}, t^8+t^9)$} &
\stack{\strut \\ $a \in \{0,1\}$ \\ modulus $64a^2-63a-56b$} \\
\hline
\end{tabular}
\end{table}

  The idea of the classification is as follows.
We first prove that parametrizations with semigroups (lexicographically) greater or equal to 
$\langle 6,7,8\rangle$ or $\langle 5,7,9\rangle$,
$\langle 5,6,13\rangle$ or
$\langle 4,9,14\rangle$  are at least bimodal  (Lemma \ref{1.3.1}). Parametrizations with smaller semigroup are  then discussed case by case and turn out to be simple (Table \ref{tab:GH}) or unimodal (Table \ref{tab:2}) .\\


\vskip 0,65 true cm
\section{Semigroups and Deformations}\label{sec2}
We start by collecting some useful properties of numerical semigroups. Recall that the generators are increasingly ordered.

For the classification of parametrizations we need  several results about semigroups. For this we order semigroups lexicographically, after making the set of generators the same length.
\begin{Definition}\label{1.2.1}
  Let $\Gamma=\langle  a_1, ..., a_l \rangle $ , $\bar{\Gamma}=\langle  b_1, ..., b_s \rangle $ be semigroups given by minimal sets of generators.
    If $l<s$ (resp. $l>s$) then we extend the set of generators $a_1, ..., a_l$ to $a_1, ..., a_l,$ $\infty, ..., \infty$ (resp. $b_1, ..., b_s$ to $b_1, ..., b_s,\infty, ..., \infty$).
Then we define $\Gamma < \bar{\Gamma}$ if $\Gamma \neq \bar{\Gamma}$
 and if there exists an $i$ such that $a_j=b_j$ for $j<i$ and
  $a_i< b_i.$\\
\end{Definition}
\begin{Example}
$\langle g_1,g_2,g_3\rangle < \langle g_1,g_2\rangle$ if $g_3\neq \infty$.
\end{Example}
The following lemma is proved in \cite[Lemma 2]{JP}. 
\begin{Lemma}\label{1.2.2}
  Let $\bar{f}=(\bar{x}(t),\bar{y}(t),\bar{z}(t))$ and ${f}=({x}(t),{y}(t),{z}(t))$ be parametrizations with semigroup $\bar{\Gamma}$ resp. $\Gamma$. If  $\Gamma \subseteq \bar {\Gamma}$ then $\bar{\Gamma} \leq {\Gamma}.$
\end{Lemma}

\begin{Lemma}\label{1.2.3}
  Let $\Gamma$ be the semigroup of the parametrization $(x(t), y(t), z(t))$ and assume that $\Gamma < \langle  5,7,16\rangle.$ 
  Let $(X(u,t), Y(u,t), Z(u,t))$ be a deformation of $(x(t), y(t), z(t))$ over $B_A$ of order $5$ and  let ${\Gamma}_{\frak{u}}$ be the corresponding semigroup for $u\in B_A$. Then there exists an open neighborhood $U$ of $0$ such that for all $u\in U$ we have ${\Gamma}_{u}\leq \Gamma$. If the order of ${\Gamma}_{\frak{u}}$ is $4$ then 
  ${\Gamma}_{u} < \langle 4,13,14\rangle $.
  \end{Lemma}
  \begin{proof}
  Let $(t^5,y(t),z(t))$ define a space curve in short parametrization with $\text{ ord}_ty(t)\leq 7$.
  Let $(X(u,t),Y(u,t),Z(u,t))$ be a deformation of $(t^5, y(t), z(t))$ over $B_A$ and  let ${\Gamma}_{\frak{u}}$ be the corresponding semigroup for $u\in B_A$. 
  If for $u\in B_A$ the order of the deformation is $5$ then 
  \begin{center}
  $\text{ ord}_tX(u,t)=5$, 
  $\text{ ord}_tY(u,t)\leq \text{ ord}_ty(t)\leq 7$ and $\text{ ord}_tZ(u,t)\leq \text{ ord}_tz(t)$.
  \end{center}
  We may assume that $X(u,t)=t^5+\sum_{i\geq 6} a(u)t^i$ with $a_i(0)=0$. Let $Y(u,t)=\sum_{i\geq 5} b_i(u)t^i$ and
 $Z(u,t)=\sum_{i\geq 5} c_i(u)t^i$ with $b_5(0)=c_5(0)=0$. Define 
 $$\bar Y(u,t):=Y(u,t)-b_5(u)X(u,t)=\sum_{i\geq 6}(b_i(u)-b_5(u)a_i(u)t^i=:\sum_{i\geq 6}\bar b_i(u)t^i.$$
 $$\bar Z(u,t):=Z(u,t)-c_5(u)X(u,t)=\sum_{i\geq 6}(b_i(u)-c_5(u)a_i(u)t^i=:\sum_{i\geq 6}\bar c_i(u)t^i.$$
 If $\bar b_6(u)\neq 0$ or $\bar c_6(u)\neq 0$ we have $5,6\in\Gamma_u$ and we are done if $\text{ ord}_ty(t)= 7$.
 If $\text{ ord}_ty(t)= 6$ then $\bar b_6(0)\neq 0$ and we may assume that $\bar b_6(u)=1$.
 Now we consider
$$\tilde Z(u,t):=\bar Z(u,t)-\bar c_6(u)\bar Y(u,t)=\sum_{i\geq 7}(\bar c_i(u)-\bar c_6(u)\bar b_i(u))t^i=:\sum_{i\geq 7}\tilde c_i(u)t^i.$$
Since $Z(0,t)=z(t)$ there exists $i$ such that $\tilde c_i(u)\neq 0$ for some $i\leq \text{ ord}_tz(t)$.
This implies that ${\Gamma}_{\frak{u}}\leq \Gamma$.\\ 
It remains to consider the case $\bar b_6(u)=\bar c_6(u)=0$. This implies that $\text{ ord}_ty(t)= 7$ and therefore $\bar b_7(0)\neq 0$. As before we may assume that $\bar b_7(u)=1$ and consider
$$\tilde Z(u,t):=\bar Z(u,t)-\bar c_7(u)\bar Y(u,t)=\sum_{i\geq 8}(\bar c_i(u)-\bar c_7(u)\bar b_i(u))t^i=:\sum_{i\geq 8}\tilde c_i(u)t^i.$$
Since $Z(0,t)=z(t)$ there exists $i$ such that $\tilde c_i(u)\neq 0$ for some $i\leq \text{ ord}_tz(t)$.
This implies that ${\Gamma}_{\frak{u}}\leq \Gamma$.\\

  If for $u\in B_A$ the order of the deformation is $\leq 4$ then we have two possibilities
  \begin{enumerate}
  \item $\text{ ord}_tX(u,t)=5$, 
  \item $\text{ ord}_tX(u,t) \leq 4$.
  \end{enumerate}
  In the first case we have $a,5 \in {\Gamma}_{u}$ for some $a\leq 4$. This implies that ${\Gamma}_{u} < \langle 4,13,14\rangle $.\\
  In the second case, if $\text{ ord}_tX(u,t) < 4$, then the claim is obviously true. Otherwise we have
  $4,a \in {\Gamma}_{u}$ with $a=\text{ ord}_tY(u,t)\leq \text{ ord}_ty(t)\leq 8$. If  
  $4$ does not devide $a$ then obviously ${\Gamma}_{u} < \langle 4,9\rangle < \langle 4,13,14\rangle$.\\
  If $a=4$ let 
  $$X(u,t)=\sum_{i\geq 4}a_i(u)t^i \text{ and } Y(u,t)=\sum_{i\geq 4}b_i(u)t^i $$
   with $a_4(0)=0$, $b_4(0)=0$, $a_5(0)\neq 0$ and $b_i(0,t)\neq 0$ for some $i\leq 8$.\\
   We may assume that $a_4(u)=u^{\alpha}\cdot e$ and $b_4(u)=u^{\beta}\cdot e'$, with $e, e'$ units. By multiplication with $e^{-1}, e'^{-1}$ we may assume that $e, e'$ are both $1$.
   
   If $\alpha\geq \beta$ we consider
   $$X(u,t)-u^{\alpha-\beta}Y(u,t)=\sum_{i\geq 5}(a_i(u)-u^{\alpha-\beta}b_i(u))t^i.$$
   This implies $4, 5\in \Gamma_u$ since $a_5(0)\neq 0$ and $b_5(0)=0$ implies $a_5(u)-u^{\alpha-\beta}b_5(u)\neq 0$.
   As before we obtain ${\Gamma}_{u} < \langle 4,13,14\rangle$.
   
   Now assume that $\alpha < \beta$ and consider
   $$\bar Y(u,t):=u^{\beta -\alpha}X(u,t)-Y(u,t)=\sum_{i\geq 5}(u^{\beta -\alpha}a_i(u)-b_i(u))t^i=:\sum_{i\geq 5} c_it^i.$$ 
   If $c_i\neq 0$ for some $i\in \{5,6,7\}$, we choose the minimal $i$ with this property. Then $4,i\in \Gamma_u$ and as before ${\Gamma}_{u} < \langle 4,13,14\rangle$.
   
   It remains to consider the case $c_5=c_6=c_7=0$ and $c_8\neq 0$. Since $\bar Y(0,t)=y(t)$ we may assume that
   $c_8=1$. Now we consider
   $$\tilde Y(u,t):=u^{2\alpha}\bar Y(u,t)-X(u,t)^2=u^{2\alpha}\sum_{i\geq 9}c_i(u)t^i-2u^{\alpha}\sum_{i\geq 5}a_i(u)t^{4+i}-(\sum_{i\geq 5}a_i(u)t^i)^2.$$ 
 Setting $u=0$ we obtain $\tilde Y(0,t)=-t^{10}$ because $a_5(0) =1$.
It follows that $9$ or $10$ is in $\Gamma_u$ and again we obtain ${\Gamma}_{u} < \langle 4,13,14\rangle$.
  \end{proof}

\begin{Lemma}\label{1.2.4}
  Let $\Gamma$ be the semigroup of the parametrization $(x(t), y(t), z(t))$ and assume that $\Gamma < \langle  4,13,14\rangle.$ 
  Let $(X(u,t), Y(u,t), Z(u,t))$ be a deformation of $(x(t), y(t), z(t))$ over $B_A$  and  let ${\Gamma}_{\frak{u}}$ be the corresponding semigroup for $u\in B_A$. Then there exists an open neighborhood $U$ of $0$ such that for all $u\in U$ we have ${\Gamma}_{u}\leq \Gamma$. 
  \end{Lemma}
    \begin{proof}
    The claim follows if $\ord_t x(t) <4$. If $\ord_t x(t) =4$, then
 we may assume that $(t^4,y(t),z(t))$ defines a short parametrization with $\text{ ord}_ty(t)\leq 11$.
  Let $(X(u,t),Y(u,t),Z(u,t))$ be a deformation of $(t^4, y(t), z(t))$ over $B_A$ and  let ${\Gamma}_{\frak{u}}$, $u\in B_A$, be the corresponding semigroup.
  If for $u\in B_A$ the order of the deformation is $\leq 3$ then obviously ${\Gamma}_{u}\leq \Gamma$.
  Assume now that the order of the deformation is $4$. We have
  \begin{center}
  $\text{ ord}_tX(u,t)=4$, 
  $a:=\text{ ord}_tY(u,t)\leq \text{ ord}_ty(t)\leq 11$, especially $4,a\in \Gamma_u$.
  \end{center}
  If $4$ does not divide $a$ then $\Gamma_u\leq\Gamma.$ 
  The cases $a=4$ and $a=8$ can be treated as in the proof of the previous lemma.
  \end{proof}

\section{Minimal Non-simple and Non-unimodal Curves}
\label{sec3}
Note that since $\K$ is an algebraically closed real valued field of characteristic $0$,
the algebraic closure $\overline \Q$ of $\Q$  is contained in $\K$ and hence $\overline \Q\{t\} \subseteq \K\{t\}$.

Let us first exclude parametrizations that are at least bimodal. 
We use {\sc Singular} to prove that the modality 
for these parametrizations
is $\ge 2$ and explain the procedure in the first case.
The  computations in  {\sc Singular} use coefficients in a finite extension of $\Q$.

\begin{Proposition}\label{1.3.1} 
Parametrizations in $t\cdot \K\{t\}^3$ of the form
 $$F=(t^l, \ t^m+{\Sigma}_{i>m,i \notin \Gamma}\ {b}_it^i, \ t^n + {\Sigma}_{i>n,i \notin \Gamma} \ {c}_it^i),\ b_i, c_i \in \K,$$
with  $\Gamma= \langle l,m,n \rangle$,  are at least bimodal in the following cases:
\begin{enumerate}
\item $(l,m,n) = (6,7,8)$,
\item  $(l,m,n) = (5,7,9)$,
\item  $(l,m,n) = (5,6,13)$,
\item  $(l,m,n) = (4,9,14)$.
\end{enumerate} 
In particular, the parametrizations   $(t^6,t^7,t^8)$, 
 $(t^5,t^7,t^{9})$, $(t^5,t^6,t^{13})$
$(t^4,t^{9},t^{14})$ and their specializations $(t^6,t^7,0)$,
$(t^5,t^7,0)$,  $(t^5,t^6,0)$, $(t^4,t^{9},0)$ are at least bimodal.
\end{Proposition}
We have always $\Gamma \subset \Gamma(F)$, but $\Gamma(F)$ may be bigger for certain values  $b_i, c_i$. For example $F=(t^6,t^7,t^8+t^9)$ has semigroup $\langle 6,7,8,17 \rangle$ 
($z^2-(yz-xz+ y^2) = 2t^{17} + t^{18}$) 
and $F=(t^4,t^9+t^{10},t^{14})$ has semigroup $\langle 4,9,14,19 \rangle$ 
($y^2-xz= 2t^{19} + t^{20}$), while for $(l,m,n) = (5,7,16)$, $F$ has alway semigroup $\langle 5,7,16 \rangle$.  We do not need these facts. 

\begin{proof} 
Note that $(t^l,t^m,0)$ is a specialization of  $(t^l,t^m,t^n)$
because $(t^l,t^m,t^n)\sim (t^l,t^m,\alpha t^n)$ for $\alpha\ne 0$.
The last statement follows, since the modality can only increase for specializations by semicontinuity of the modality (Proposition \ref{prop.sem}).\\
Since  the condition of Proposition \ref{prop.Za} (2) is fulfilled for $ (l,m,n) = (6,7,8), (5,7,13)$, and $(4,9,14)$, we have to show that
the parametrization 
$$F=(t^l, \ t^m+{\Sigma}_{c-1>i>m,i \notin \Gamma}\ {b}_it^i, \ t^n + {\Sigma}_{c-1>i>n,i \notin \Gamma} \ {c}_it^i),$$
with arbitrary $b_i$ and $c_i$ in $\K$, is at least bimodal. We will also use Remark \ref{rm.Za}, saying that we can make the second non-vanishing term in the second or third component of a parametrization equal to 1.\\
For $ (l,m,n) = (5,7,9)$ the condition of Proposition \ref{prop.Za} (2) is not fulfilled  for the  second component and thus we have to keep $b_{c-1}=b_{13}$.

1. Let  
$\Gamma = \langle 6,7,8 \rangle$ with gaps  $1,2,3,4,5,9,10,11,17$ and $c=18$. The short parametrization is (using Remark \ref{rm.Za})
$$F=F(a,b,t)=(t^6,\ t^7+t^9+at^{10}+bt^{11},\ t^8+et^{10}+ft^{11})$$ 
with $a,b,e,f  \in \K$. 
We will use {\sc Singular} below to prove that  for fixed  $e,f\in \mathbb K$.
$$F(a,b,t) \sim F(c,d,t)$$
implies  $(a,b)=(c,d)$. The computations are in the ring $S$ below, with $e,f$ generic.
This shows that the family $F(a,b,t)$,  $(a,b) \in  \K^2$, is modular. By Proposition \ref{prop.modular} 
this implies that the parametrization $F$, with  $e,f$ generic, is at least bimodal. Since $F$ with non-generic elements $e,f \in \K$ is a specialization, it follows from the semicontinuity of the modality that 
$F$ is  at least bimodal for arbitrary  $a,b,ef, \in \K$.\\

Let 
$$
\begin{array}{lll}
(x(t,y(t),z(t)) )&=&(t^6, \ t^7+t^9+at^{10}+bt^{11}, \ t^8+et^{10}+ft^{11})  \text{ and }\\
 (\hat x(t,\hat y(t),\hat z(t)) )&=&(t^6, \ t^7+t^9+ct^{10}+dt^{11}, \ t^8+et^{10}+ft^{11})
 \end{array}
$$
Using Remark \ref{rm.aut} we may consider as Ansatz
 the equations 
 \begin{align}\tag{*}\label{*}
 \begin{split}
x(\psi(t))&=  H_x(\hat{x}(\alpha t),\hat{y}(\alpha t),\hat{z}(\alpha t)),\\
y(\psi(t))&= H_y(\hat{x}(\alpha t),\hat{y}(\alpha t),\hat{z}(\alpha t)),\\
z(\psi(t))&= H_z(\hat{x}(\alpha t),\hat{y}(\alpha t),\hat{z}(\alpha t)),
\end{split}
\end{align}
with $\psi(t)=t+a_2t^2+a_3t^3+\ldots$ and $H_x,H_y,H_z$ defining an automorphism of $\mathbb K\{x,y,z\}$.
We consider these equations modulo $t^{11}\K\{t\}$.

The {\sc Singular} code reads as follows:
\begin{verbatim}
ring R=(0,alpha,a2,a3,a4,a5,a6,a7,a8,a9,a10,a11,a12,a13,a14,a15,a16,u1,u2,
u3,u4,u5,u6,u7,u8,u9,u10,u11,u12,u13,u14,v1,v2,v3,v4,v5,v6,v7,v8,v9,v10,
v11,v12,v13,v14,w1,w2,w3,w4,w5,w6,w7,w8,w9,w10,w11,w12,w13,w14,z1,z2,z3,
z4,z5,z6,z7,z8,z9,z10,z11,z12,z13,z14,a,b,c,d,e,f),(x,y,z,t),ds;
poly psi=t+a2*t2+a3*t3+a4*t4+a5*t5+a6*t6+a7*t7+a8*t8+a9*t9+a10*t10+a11*t^11
         +a12*t^12+a13*t^13+a14*t^14+a15*t^15+a16*t^16;
poly psi1=t+a2*t2+a3*t3+a4*t4+a5*t5+a6*t6+a7*t7+a8*t8+a9*t9+a10*t10+a11*t11;      
poly psi2=t+a2*t2+a3*t3+a4*t4+a5*t5+a6*t6;
poly Hx=u1*x+u2*y+u3*z+u5*x2+u6*xy+u7*xz+u9*y2+u10*yz+u12*z2+u13*x3+
        u14*x2y+u4*x^4;
poly Hy=v1*x+v2*y+v3*z+v5*x2+v6*xy+v7*xz+v9*y2+v10*yz+v12*z2+v13*x3+
        v14*x2y+v4*x^4;
poly Hz=w1*x+w2*y+w3*z+w5*x2+w6*xy+w7*xz+w9*y2+w10*yz+w12*z2+w13*x3+
        w14*x2y+w4*x^4;

poly X=jet(psi^6
-subst(Hx,x,(alpha*t)^6,y,(alpha*t)^7+(alpha*t)^9+a*(alpha*t)^10
      +b*(alpha*t)^11,z,(alpha*t)^8+e*(alpha*t)^10+f*(alpha*t)^11),11);  
poly Y=jet(psi^7+psi^9+c*psi^10+d*psi^11
-subst(Hy,x,(alpha*t)^6,y,(alpha*t)^7+(alpha*t)^9+a*(alpha*t)^10
      +b*(alpha*t)^11,z,(alpha*t)^8+e*(alpha*t)^10+f*(alpha*t)^11),11);  
poly Z=jet(psi^8+e*psi^10+f*psi^11
-subst(Hz,x,(alpha*t)^6,y,(alpha*t)^7+(alpha*t)^9+a*(alpha*t)^10
      +b*(alpha*t)^11,z,(alpha*t)^8+e*(alpha*t)^10+f*(alpha*t)^11),11);  
matrix M1=coef(X,t);matrix M2=coef(Y,t);matrix M3=coef(Z,t);
ideal I;int ii;
for(ii=1;ii<=ncols(M1);ii++){I[size(I)+1]=M1[2,ii];}
for(ii=1;ii<=ncols(M2);ii++){I[size(I)+1]=M2[2,ii];}
for(ii=1;ii<=ncols(M3);ii++){I[size(I)+1]=M3[2,ii];}
\end{verbatim}
We created the ring $R$, which is the localization (because of ordering {\tt ds}) at the ideal $\langle x,y,z,t\rangle$ of the polynomial ring with variables $x,y,z,t$ and parameters 
$\alpha,a_2,a_3,...,m,n$, i.e. the localization
$\Q(\alpha,a_2,a_3,...,m,n)[x,y,z,t]_{\langle x,y,z,t\rangle}$. The Ansatz for $\psi$ with parameters $a_2,...,a_{10}$ and for $H_x,H_y,H_z$ with parameters $u_i,v_i,w_i, i=1,...,10$ provides an $\mathcal A$-equivalence up to order 10, iff
the system (\ref{*}) has a solution up to order 10, i.e., if the polynomial system $X=0, Y=0, Z=0$ can be solved. This is equivalent that the system of equations given by the $t$-coefficients, which are collected in the ideal $I$, can be solved. \\
To check this, we map $I$ to the  polynomial ring $S$ and compute a Gr\"obner basis:
\begin{verbatim}
ring S=(0,e,f),(a2,a3,a4,a5,a6,a7,a8,a9,a10,a11,a12,a13,a14,a15,a16,u1,
u2,u3,u4,u5,u6,u7,u8,u9,u10,u11,u12,u13,u14,v1,v2,v3,v4,v5,v6,v7,v8,v9,v10,
v11,v12,v13,v14,w1,w2,w3,w4,w5,w6,w7,w8,w9,w10,w11,w12,w13,w14,z1,z2,z3,z4,
z5,z6,z7,z8,z9,z10,z11,z12,z13,z14,a,b,c,d,alpha),lp;
ideal I=imap(R,I);
option(redSB);
ideal J=std(I);

J[1]=alpha-1
J[2]=b-d
J[3]=a-c
\end{verbatim}
It follows that $F(a,b,t) \sim F(c,d,t)$
implies $a=c$ and $b=d$.
This shows that the 10-jet of the family $F$ and hence $F$ itself is modular, showing that every parametrization with semigroup $\Gamma=\langle 6,7,8 \rangle$ is at least bimodal, as claimed. \\

We use the rings $R$ and $S$ also for the following computations.\\

2.  Let  
$\Gamma = \langle 5,7,9 \rangle$ with gaps 
$1,2,3,4,6,8,11,13$ and $c=14$. Similar to the case before we use {\sc{Singular}} to prove that  the equivalence of the short parametrizations
$$F(a,b,t)=(t^5,t^{7}+at^{8}+bt^{11}+et^{13},t^{9}+ft^{11})\sim (t^5,t^{7}+ct^{8}+dt^{11}+et^{13},t^{9}+ft^{11})
$$
with fixed $e,f  \in \K$ 
implies $a=c$ or $a=-c$ and $b=d$:

\begin{verbatim}
setring R;
X=jet(psi^5
-subst(Hx,x,(alpha*t)^5,y,(alpha*t)^5,y,(alpha*t)^7+a*(alpha*t)^8
      +b*(alpha*t)^11+e*(alpha*t)^13,z,(alpha*t)^9+f*(alpha*t)^11),13);
Y=jet(psi^7+c*psi^8+d*psi2^11+e*psi2^13
-subst(Hy,x,(alpha*t)^5,y, (alpha*t)^7+a*(alpha*t)^8+b*(alpha*t)^11
      +e*(alpha*t)^13,z,(alpha*t)^9+f*(alpha*t)^11),13);
Z=jet(psi^9+f*psi2^11
-subst(Hz,x,(alpha*t)^5,y,(alpha*t)^7+a*(alpha*t)^8+b*(alpha*t)^11
      +e*(alpha*t)^13,z,(alpha*t)^9+f*(alpha*t)^11),13);
M1=coef(X,t);  M2=coef(Y,t);  M3=coef(Z,t);
ideal I;
for(ii=1;ii<=ncols(M1);ii++){I[size(I)+1]=M1[2,ii];}
for(ii=1;ii<=ncols(M2);ii++){I[size(I)+1]=M2[2,ii];}
for(ii=1;ii<=ncols(M3);ii++){I[size(I)+1]=M3[2,ii];}
setring S;
ideal I=imap(R,I);
ideal J=std(I);

J[1]=alpha^2-1
J[2]=b-d
J[3]=a-c*alpha
\end{verbatim}
This shows again that the 2-dimensional family $F(a,b,t)$ is modular and hence that parametrizations with semigroup $\langle 5,7,9 \rangle$ are at least bimodal. \\

3.  Let  
$\Gamma = \langle 5,6,13 \rangle$ with gaps 
$1,2,3,4,7,8,9,14$ and $c=15$. Similar to the case before we use {\sc{Singular}} to prove that  the equivalence of the short parametrizations
$$(t^5,t^6+t^{7}+at^{8}+bt^{9},t^{13})\sim (t^5,t^6+t^{7}+ct^{8}+dt^{9},t^{13})
$$
implies $a=c$ and $b=d$:

\begin{verbatim}
setring R;
X=jet(psi^5
-subst(Hx,x,(alpha*t)^5,y,(alpha*t)^6+(alpha*t)^7
      +a*(alpha*t)^8+b*(alpha*t)^9,z,(alpha*t)^13),14); 
Y=jet(psi^6+psi^7+c*psi^8+d*psi^9
-subst(Hy,x,(alpha*t)^5,y,(alpha*t)^6(alpha*t)^7
      +a*(alpha*t)^8+b*(alpha*t)^9,z,(alpha*t)^13),14); 
Z=jet(psi2^13
-subst(Hz,x,(alpha*t)^5,y,(alpha*t)^6+(alpha*t)^7
      +a*(alpha*t)^8+b*(alpha*t)^9,z,(alpha*t)^13),14); 
M1=coef(X,t);  M2=coef(Y,t);  M3=coef(Z,t);
ideal I;
for(ii=1;ii<=ncols(M1);ii++){I[size(I)+1]=M1[2,ii];}
for(ii=1;ii<=ncols(M2);ii++){I[size(I)+1]=M2[2,ii];}
for(ii=1;ii<=ncols(M3);ii++){I[size(I)+1]=M3[2,ii];}
setring S;
ideal I=imap(R,I);
ideal J=std(I);

J[1]=alpha-1
J[2]=b-d
J[3]=a-c
\end{verbatim}

4. Let  
$\Gamma = \langle 4,9,14 \rangle$ with gaps 
$1,2,3,5,6,7,10,11,15,19$ and $c=20$. Similar to the previous cases we use {\sc{Singular}} to prove that  for each $f\in \K$ 
$$(t^4,t^{9}+t^{10}+at^{11}+ft^{15},t^{14} +bt^{15})\sim (t^4,t^{9}+t^{10}+ct^{11}+ft^{15},t^{14} +dt^{15})
$$
implies $a=c$ and $b=d$.
\begin{verbatim}
setring R;
X=jet(psi^4
-subst(Hx,x,(alpha*t)^4,y,(alpha*t)^9+(alpha*t)^10
      +a*(alpha*t)^11+ f*(alpha*t)^15,z,(alpha*t)^14+b*(alpha*t)^15),19); 
Y=jet(psi1^9+psi1^10+c*psi1^11+f*psi2^15
-subst(Hy,x,(alpha*t)^4,y,(alpha*t)^9+(alpha*t)^10
      +a*(alpha*t)^11+f*(alpha*t)^15,z,(alpha*t)^14+b*(alpha*t)^15),19); 
Z=jet(psi2^14+d*psi2^15
-subst(Hz,x,(alpha*t)^4,y,(alpha*t)^9+(alpha*t)^10
      +a*(alpha*t)^11+f*(alpha*t)^15,z,(alpha*t)^14+b*(alpha*t)^15),19);  
M1=coef(X,t);M2=coef(Y,t);M3=coef(Z,t);
ideal I;
for(ii=1;ii<=ncols(M1);ii++){I[size(I)+1]=M1[2,ii];}
for(ii=1;ii<=ncols(M2);ii++){I[size(I)+1]=M2[2,ii];}
for(ii=1;ii<=ncols(M3);ii++){I[size(I)+1]=M3[2,ii];}
setring S;
I=imap(R,I); 
ideal J=std(I);

J[1]=alpha-1
J[2]=d-b
J[3]=a-c
\end{verbatim}
\end{proof}

\begin{Proposition}\label{1.3.2} 
Let $F, G \in t\cdot \K\{t\}^3$ be parametrizations of the form
 $$
 \begin{array}{llll}
 F&=&(t^l, \ t^m+{\Sigma}_{i>m,i \notin \Gamma}\ {b}_it^i, \ t^n + {\Sigma}_{i>n,i \notin \Gamma} \ {c}_it^i),\ b_i, c_i \in \K,
  &\Gamma= \langle l,m,n \rangle,\\
 G&=&(t^{l'}, \ t^{m'}+{\Sigma}_{i>m',i \notin \Gamma'}\ {b'}_it^i, \ t^{n'} + {\Sigma}_{i>n',i \notin \Gamma'} \ {c'}_it^i),\ b'_i, c'_i \in \K,
  &\Gamma'= \langle l',m',n' \rangle,
  \end{array}
   $$
 Assume that $l\le l',m\le m',n\le n'$, and that  the $\mathcal A$-modality of $F$ is bigger or equal to $M$ for all choices of $b_i,c_i$.
 Then the $\mathcal A$-modality of $G$ is bigger or equal to $M$ for all choices of $b'_i,c'_i$.
 \end{Proposition}
 \begin{proof} 
 Write  $F=(x_1(t),x_2(t),x_3(t))$ resp. $G=(y_1(t),y_2(t),y_3(t))$ and let  $h_i \in \{l,m,n\}$ resp.  $k_i \in \{l',m',n'\}$ be the  leading exponent of $x_i(t)$ resp. $y_i(t)$. 
 Define 
 $$F'=(x'_1(t),x'_2(,t),x'_3(t))=(a_1t^l +t^{l'}, \ a_2t^m+y_2(t), \ a_3t^n + y_3(t))$$
with $a_i=0$ if $h_i=k_i$ and  $a_i=1$ if $h_i<k_i$. 
By Theorem \ref{thm.Za} $F'$ is $\mathcal A$-equivalent to 
$F$ (for some coefficients $b_i,c_i$) and hence has modality $\ge M$. Moreover, choosing $a_i$ arbirary in $\K\setminus \{ 0\}$, we get a parametrization $G'$, which  is equivalent to $F'$. We see that  $G$ ($a_i=0$) is a specialization of $G'$ and
by semicontinuity (Proposition \ref{prop.sem}) the modality of $G$  is $\ge M$.
 \end {proof}

\section{Curves of order 4} \label{sec4}
In this and the next section we classify parametrizations of order 4 and 5 by providing normal forms (which are of course not unique).
Let $f=(x(t),y(t),z(t))$ define a unimodal parametrized curve singularity and $\Gamma$ its semigroup. 
By Corollary  \ref{cor.cho} the order of $f$ is at least $4$ and 
$\Gamma > \langle 4,7\rangle$.  

If the order of $f$ is $4$ and $\Gamma > \langle 4,7\rangle$, then 
$\langle 4,a\rangle \subset \Gamma$, with $a\ge 9$ a gap of $\langle 4,7\rangle$.  
We obtain from Proposition \ref{1.3.1}  that $\Gamma <\langle 4,9,14\rangle$. 
Gaps $> 9$ of $\langle 4,9\rangle$ are 
10, 11, 14, 15, 19, 23.
Since 
$\langle 4,9,10,11\rangle$ and
$\langle 4,9,10,15\rangle$ cannot be realized as space curves,
we have to classify only curves with semigroup $\langle 4,9,10\rangle$ and $\langle 4,9,11\rangle$.

\begin{Lemma}\label{lem.4910}
Let $f=(x(t),y(t),z(t))$ be a parametrization with semigroup $\langle 4,9,10\rangle$. Then $f$ is equivalent to one of  the following normal forms:
\begin{itemize}
\item $(t^4,t^9+t^{11},t^{10}+at^{11})$, $a\in \K$,
\item $(t^4,t^9,t^{10}+at^{11})$ with $a\in \{0,1\}$.
\end{itemize}
\end{Lemma}
\begin{proof} The semigroup $\langle 4,9,10\rangle$ has gaps $1,2,3,5,6,7,11,15$.
The short parametrization (Proposition \ref{prop.Za}) shows
$$(x(t),y(t),z(t))=(t^4,t^9+at^{11}+bt^{15},t^{10}+ct^{11}+dt^{15})
\sim (t^4,t^9+at^{11},t^{10}+ct^{11})$$ 
and by and Remark \ref{rm.Za} we get
$$(t^4,t^9+at^{11},t^{10}+ct^{11})\sim (t^4,t^9+t^{11},t^{10}+ct^{11}), c\in \K \text{ if } a\ne 0$$ or 
$\sim (t^4,t^9,t^{10}+ct^{11})$, $c\in \{0,1\}$ if $a=0$.
\end{proof}
An analogous argument for  the semigroup $\langle 4,9,11\rangle$ (gaps $1,2,3,5,6,7,10,14$) gives the following normal forms.
\begin{Lemma}\label{lem.4911}
Let $f=(x(t),y(t),z(t))$ be a parametrization with semigroup $\langle 4,9,11\rangle$. Then $f$ is equivalent to one of  the following normal forms:
\begin{itemize}
\item $(t^4,t^9+t^{10},t^{11}+at^{14})$, $a\in \K$,
\item $(t^4,t^9,t^{11}+at^{14})$ with $a\in \{0,1\}$.
\end{itemize}
\end{Lemma}

\section{Curves of order 5} \label{sec5}
Let $f=(x(t),y(t),z(t))$ be a unimodal parametrization of order $\ge 5$.  Using Corollary  \ref{cor.cho} we have to consider
$\Gamma \ge \langle 4,7\rangle$. In Section \ref{sec4} we treated order 4.
If the order of $f$ is $\ge 6$, it is a specialization of the parametrization $(t^6,t^7,t^8)$, which is at least bimodal by
Proposition \ref{1.3.1}. By semicontinuity of the modality $f$ itself is at least bimodal. 

Let the order of $f$ be 5. Proposition \ref{1.3.1}
 implies that a parametrization $f$ with semigroup $\Gamma(f) = \langle 5,6,13\rangle$ is at least bimodal. If $\Gamma(f) > \langle 5,6,13\rangle$ then $f=(t^5,y(t),z(t))$ is a specialization of 
 $(t^5,at^6+y(t),bt^{13}+z(t))$ ($a\ne0$ if $\ord(y)>6$ and $b\ne0$ if $\ord(z)>13$). Since the latter is at least bimodal, $f$ is at least bimodal. The same holds if $\Gamma(f)\geq \langle 5,7,9\rangle$.\\
We have thus to consider parametrizations $f$ of order $5$ with $\Gamma (f)< \langle 5,6,13\rangle$ and $< \langle 5,7,9\rangle$. Gaps $>5$ of $\langle 5,6,13\rangle$ are $7,8,9,14$ and Gaps $>5$ of $\langle 5,7,8\rangle$ are $6,9,11$.
 
\begin{Lemma} \label{lem.5716}
Let  $\Gamma < \langle 5,7,13\rangle$ be a semigroup of order 5, such that there exists a parametrization  $f=(x(t),y(t),z(t))$ of a space curve with $\Gamma = \Gamma(f)$.  If $f$ is at most unimodal then  $\Gamma$ is one of the following semigroups.
\begin{enumerate}
\item $\langle 5,6\rangle \subset \Gamma$. Then
$\Gamma \in \{\langle 5,6\rangle, \langle 5,6,7\rangle,
 \langle 5,6,8\rangle,  \langle 5,6,9\rangle\}$. 
\item $\langle 5,7\rangle \subset \Gamma$. Then
$\Gamma \in \{\langle 5,7\rangle, \langle 5,7,8\rangle\}$.
\end{enumerate}
\end{Lemma}
\begin{proof} 1. 
We may assume that
$$f=(x(t),y(t),z(t))=(t^5,\ t^6+a_1t^{8}+a_2t^{9}+..., \ t^{n}+...)$$
with $n\ge 7$ and $n \notin \langle 5,6\rangle$ (Proposition \ref{prop.Za}). 
The gaps $\ge 7$ of $\langle 5,6\rangle$ are $n=7, 8, 9, 13, 14, 19$
and  $\Gamma = \langle 5,6,n\rangle$ has conductor
 $c=10,10,14,15,14,15$ respectively. 
Setting $(x,y,z)=(x_1,x_2,x_3)$, with orders $n_1=5,n_2 =6,n_3 =n$, we see that the order of $x_ix_j$ is either bigger than $c$ or different from the order of another product $x_kx_l$. That is, there is no cancellation among such products and we cannot produce new generators of $\Gamma(f)$ other than 
 $\langle 5,6,n\rangle $.\\
2. Let 
$$f=(x(t),y(t),z(t))=(t^5,\ t^7+a_1t^{8}+a_2t^{9}+..., \ t^{n}+...)$$ with $n\ge 8$ and $n \notin \langle 5,7\rangle$. The gaps $\ge 8$ of $\langle 5,7\rangle$ are $n=8, 9, 11, 13, 16, 18, 23$
and  the semigroup $\Gamma = \langle 5,7,n\rangle$ has conductor
 $c=12,14,14,17,19,17,19$ respectively. 
 By Proposition \ref{1.3.1} $(t^5,t^7,t^{9})$ has modality $\ge 2$ and it follows from Proposition \ref{1.3.2} that $f$ with $n\ge 9$ has also modality $\ge 2$.
It follows that $\Gamma(f) =\langle 5,7\rangle $ or $\langle 5,7,8\rangle $.
\end{proof}

By Lemma \ref{lem.5716} we have to classify parametrization $f$ of order 5 with 
 $\langle 5,6\rangle, \langle 5,7\rangle \subset \Gamma(f).$
 
\subsection {The case  $\langle 5,6\rangle \subset \Gamma$}
We have to classify parametrization $f=(x(t),y(t),z(t))$  with
semigroups $\Gamma(f)$, where $\Gamma$ is one of the semigroups 
from $ \{\langle 5,6\rangle, \langle 5,6,7\rangle,
 \langle 5,6,8\rangle,  \langle 5,6,9\rangle\}$ (Lemma \ref{lem.5716}.1.)
\begin{enumerate}
\item [0.] The semigroup $\langle 5,6\rangle$.\\
The plane curve singularities with semigroup
$\langle 5,6\rangle$ are classified in \cite[Lemma 4]{MP}.
As plane curves they are unimodal. However, the space curve
$(t^5,t^6+ \text{ higher terms },0)$ is not unimodal, since it is a specialization of
$(t^5,t^6+ \text{ higher terms },t^{13})$, which has modality 
$\ge 2$ by Proposition \ref{1.3.1}.
\item The semigroup $\langle 5,6,7\rangle$. 

If $\Gamma(f) =  \langle 5,6,7\rangle$ then $f$ is $\mathcal A$-equivalent to one of  the following normal forms:
\begin{itemize}
\item
$(t^5,t^6+t^8+at^9,t^7)$, $a \in \K$.
\item
$(t^5,t^6+at^9,t^7) \text{ with } a\in\{0,1\}.$
\end{itemize}
The gaps $>6$ of $\langle 5,6,7\rangle$ are 8, 9 and c=10. By Proposition \ref{prop.Za} and Remark \ref{rm.Za} we may assume that
$$ (x(t),y(t),z(t)) =(t^5,t^6+ct^7+dt^{8}+et^{9},t^7+ft^{8} )$$
To deduce the normal forms we use {\sc{Singular}} to prove first that  for given $c,d,e,f\in \mathbb K$ there exist $a,b\in \mathbb K$ such that
$$(t^5,t^6+ct^7+dt^{8}+et^{9},t^7+ft^{8} )\sim (t^5,t^6+at^8+bt^9,t^7).
$$
Let $$x(t)=\hat x(t)=t^5, y(t)=t^6++ct^7+dt^8+et^9, \hat y(t)=t^6+at^{8}+bt^{9} $$ and 
$$ z(t)=t^7+ft^{8}, \hat z(t)=t^7 .
$$
Using Remark \ref{rm.aut} we may consider as Ansatz
 the equations
$$x(\psi(t))=H_x(\hat{x}(\alpha t),\hat{y}(\alpha t),\hat{z}(\alpha t)),$$
$$y(\psi(t))=H_y(\hat{x}(\alpha t),\hat{y}(\alpha t),\hat{z}(\alpha t)),$$
$$z(\psi(t))=H_z(\hat{x}(\alpha t),\hat{y}(\alpha t),\hat{z}(\alpha t)).$$
with $\psi(t)=t+a_2t^2+a_3t^3+\ldots$ and $H_x,H_y,H_z$ defining and automorphism of $\mathbb K[[x,y,z]]$.\\
Using Lemma \ref{lem.determined} it is enough to consider the equations modulo the conductor ideal.

The following {\sc Singular} computationshows that the equations for the Ansatz have a solution.
\begin{verbatim}
ring R=(0,a1,a2,a3,a4,a5,a6,a7,a8,a9,a10,u1,u2,u3,u4,u5,u6,u7,u8,u9,
     u10,u11,u12,u13,u14,v1,v2,v3,v4,v5,v6,v7,v8,v9,v10,v11,v12,v13,
     v14,w1,w2,w3,w4,w5,w6,w7,w8,w9,w10,w11,w12,w13,w14,z1,z2,z3,z4,
     z5,z6,z7,z8,z9,z10,z11,z12,z13,z14,a,b,c,d,e,f,g,h,i,j,k,l,m,n),
     (x,y,z,t),ds;
poly psi=t+a2*t2+a3*t3+a4*t4+a5*t5+a6*t6+a7*t7+a8*t8+a9*t9+a10*t10;
poly Hx=u1*x+u2*y+u3*z+u5*x2+u6*xy+u7*xz+u9*y2+u10*yz+u12*z2;
poly Hy=v1*x+v2*y+v3*z+v5*x2+v6*xy+v7*xz+v9*y2+v10*yz+v12*z2;
poly Hz=w1*x+w2*y+w3*z+w5*x2+w6*xy+w7*xz+w9*y2+w10*yz+w12*z2;
poly X=jet(psi^5
    -subst(Hx,x,(a1*t)^5,y,(a1*t)^6+c*(a1*t)^7+d*(a1*t)^8+e*(a1*t)^9,z,
    (a1*t)^7+f*(a1*t)^8),10);
poly Y=jet(psi^6+a*psi^8+b*psi^9
    -subst(Hy,x,(a1*t)^5,y,(a1*t)^6+c*(a1*t)^7+d*(a1*t)^8+e*(a1*t)^9,z,
    (a1*t)^7+f*(a1*t)^8),10);
poly Z=jet(psi^7
    -subst(Hz,x,(a1*t)^5,y,(a1*t)^6+c*(a1*t)^7+d*(a1*t)^8+e*(a1*t)^9,z,
    (a1*t)^7+f*(a1*t)^8) ,10);
matrix M1=coef(X,t);matrix M2=coef(Y,t);matrix M3=coef(Z,t);
ideal I;int ii;
for(ii=1;ii<=ncols(M1);ii++){I[size(I)+1]=M1[2,ii];}
for(ii=1;ii<=ncols(M2);ii++){I[size(I)+1]=M2[2,ii];}
for(ii=1;ii<=ncols(M3);ii++){I[size(I)+1]=M3[2,ii];}
I=subst(I,a1,1);

ring S1=(0,c,d,e,f),(a2,a3,a4,a5,a6,a7,a8,a9,a10,u1,u2,u3,u4,u5,u6,
u7,u8,u9,u10,u11,u12,u13,u14,v1,v2,v3,v4,v5,v6,v7,v8,v9,v10,v11,v12,
v13,v14,w1,w2,w3,w4,w5,w6,w7,w8,w9,w10,w11,w12,w13,w14,z1,z2,z3,z4,z5,
z6,z7,z8,z9,z10,z11,z12,z13,z14,g,h,i,j,k,l,m,n,a1,a,b),lp;
ideal I=imap(R,I);
ideal J=std(I);
dim(J);         
// 57
\end{verbatim}
The dimension of the zero set of the equations corresponding to the Ansatz is positive. Therefore we obtain
a solution for the Ansatz and hence  $f \sim (t^5,t^6,t^7+at^{8}+bt^{9})$  for generic $c,d,e,f$ (in the ring S1 we compute with polynomials over $\Q(c,d,e,f)$). We show now that indeed $f \sim (t^5,t^6,t^7+at^{8}+bt^{9})$ for all $c,d,e,f \in \Q$. 
Looking at the first two elements $F,G$ of $J$, we see that $a$ and $b $ are simple expressions in 
$\Q[c,d,e,f]$, which hold in $I\cdot \Q(c,d,e,f)[a_2 ...,u_1...,...]$.
\begin{verbatim}
poly F=J[1];
poly G=J[2];
F;
// 49*b+(-98*c*f^2+98*d*f-49*e+44*f^3)
G;
// 49*a+(49*c*f-49*d-45*f^2)
\end{verbatim}
In the next step we show that the expressions for $F,G$ hold in  in $I \subset \Q[c,d,e,f,a_2 ...,u_1...,...].$
\begin{verbatim}
ring S=0,(c,d,e,f,a2,a3,a4,a5,a6,a7,a8,a9,a10,u1,u2,u3,u4,u5,u6,u7,u8,u9,u10,
u11,u12,u13,u14,v1,v2,v3,v4,v5,v6,v7,v8,v9,v10,v11,v12,v13,v14,w1,w2,
w3,w4,w5,w6,w7,w8,w9,w10,w11,w12,w13,w14,z1,z2,z3,z4,z5,z6,z7,z8,z9,
z10,z11,z12,z13,z14,a,b,g,h,i,j,k,l,m,n),lp;
ideal I=imap(R,I);
ideal J=std(I);
poly F=imap(S1,F);
poly G=imap(S1,G);
reduce(F,J);
// 0
reduce(G,J);
//0
\end{verbatim}
Hence $f\sim (t^5,t^6+at^8+bt^9,t^7)$.
If $a=0$ we get the second normal form $(t^5,t^6+bt^{9},t^7)$, with $b\in \{0,1\}$ by Remark \ref{rm.Za}. If $a\ne 0$ we can make $a=1$ by Remark \ref{rm.Za} and we get the first normal form $(t^5,t^6+t^8+bt^{9},t^7)$.

For the remaining semigroups we obtain the following normal forms without using \sc{Singular}.\\

2. The semigroup $\Gamma = \langle 5,6,8\rangle$
\begin{itemize}
\item $(t^5,t^6+at^7,t^8)$ with  $a\in\{0,1\}$.
\end{itemize}
The gaps $>6$ are 7 and 9. The short parametrization $(t^5,t^6+at^7+bt^9,t^8 +ct^9)$ can be reduced to 
$(t^5,t^6+at^7,t^8)$ by Proposition \ref{prop.Za} and Remark \ref{rm.Za}, with $a\in \{0,1\}$. This is unimodal since it is a specialization of  $(t^5,t^6+at^7,\alpha t^7+t^8)$, which has semigroup $\langle 5,6,7\rangle$ and is unimodal by case 1.

3. The semigroup $\Gamma = \langle 5,6,9\rangle$
\begin{itemize}
\item $(t^5,t^6+t^7+bt^8, t^9)$ with $b \in \K$.
\item $(t^5,t^6+bt^8, t^9)$ with $b \in \{0,1\}$.
\end{itemize}
The gaps $>6$ are 7,8,13. The short parametrization 
$(t^5,t^6+at^7+bt^8 +ct^{13}, t^9+ dt^{13})$ can be reduced to 
$(t^5,t^6+at^7+bt^8, t^9)$ by Proposition \ref{prop.Za} and Remark \ref{rm.Za}, with $a\in \{0,1\}$ and $b\in \{0,1\}$ if $a=0$.
\end{enumerate}

\subsection {The case $\langle 5,7\rangle \subset \Gamma$}
According to Lemma \ref{lem.5716}, we have to consider only the semigroups 
$\langle 5,7\rangle$, $\langle 5,7,8\rangle$. 
\begin{enumerate}
\item [0.] The semigroup $\Gamma =\langle 5,7\rangle$\\
The plane curve singularities with semigroup
$\langle 5,7\rangle$ are classified in \cite[Lemma 4]{MP}.
As plane curves they are unimodal. However, the space curve
$(t^5,t^7+ \text{ higher terms },0)$ is not unimodal, since it is a specialization of
$(t^5,t^7+ \text{ higher terms },t^{13})$, which has modality 
$\ge 2$ by Proposition \ref{1.3.1}.

\item The semigroup $\Gamma = \langle 5,7,8\rangle$

If $\Gamma(f) =  \langle 5,7,8\rangle$ then $f$ is $\mathcal A$-equivalent to one of  the following normal forms:\\
\begin{itemize}
\item  $(t^5,t^7+bt^{11},t^8+ct^{9})$ with $c\in \{0,1\}, b\in \K$,
\item
$(t^5,t^7+at^9+bt^{11},t^{8}+t^{9})$ with modulus $64a^2-63a-56b.$\\
\end{itemize}

The gaps $>8$ of $\langle 5,7,8\rangle$ are 9, 11 and c=12. By 
Theorem \ref{thm.Za} we may assume that $f$ is equal to 
$f=(t^5,t^7+at^{9}+bt^{11},t^8+ct^{9}+dt^{11} )$
and by Proposition \ref{prop.Za} 
 $f$ is equivalent to 
$$(t^5,t^7+at^{9}+bt^{11},t^8+ct^{9}).$$
$a=0$:  By Remark \ref{rm.Za}  this is equivalent to 
$(t^5,t^7+bt^{11},t^8+ct^{9})$ with $c\in \{0,1\}, b\in \K$ arbitrary.\\
$a\ne 0$: By Remark \ref{rm.Za} this is equivalent to 
 $(t^5,t^7+at^{9}+bt^{11},t^8+t^{9})$ with $a,b\in \K$
(or to $(t^5,t^7+t^{9}+bt^{11},t^8+ct^{9})$ with $b,c\in \K$.)

Consider the family of parametrizations $F(a,b,t)=(t^5,t^7+at^9+bt^{11},t^8+t^9)$.
We use {\sc Singular} to show that
$$(t^5,t^7+at^9+bt^{11},t^8+t^9)\sim
 (t^5,t^7+mt^9+nt^{11},t^8+t^9)$$
is equivalent to $H(a,b) =H(m,n)$ with $H(u,v) =64u^2-63u-56v$, defining a curve in the $(u,v)$-plane.
So the equivalence classes in the $(u,v)$-plane are the fibers of the map 
$H:\K^2 \to \K$, which is a geometric quotient for the $\mathcal A$-action on $K^2$ (Definition \ref{def.modal}). Hence $F(a,b,t)$ is unimodal, with modulus $64a^2-63a-56b.$

By Lemma \ref{lem.determined} it is enough to consider the equations modulo the conductor ideal.
\begin{verbatim}
ring R=(0,a1,a2,a3,a4,a5,a6,a7,a8,a9,a10,u1,u2,u3,u4,u5,u6,u7,u8,u9,
u10,u11,u12,u13,u14,v1,v2,v3,v4,v5,v6,v7,v8,v9,v10,v11,v12,v13,v14,
w1,w2,w3,w4,w5,w6,w7,w8,w9,w10,w11,w12,w13,w14,z1,z2,z3,z4,z5,z6,z7,
z8,z9,z10,z11,z12,z13,z14,a,b,c,d,e,f,g,h,i,j,k,l,m,n),(x,y,z,t),ds;
poly psi=t+a2*t2+a3*t3+a4*t4+a5*t5+a6*t6+a7*t7+a8*t8+a9*t9+a10*t10;
poly Hx=u1*x+u2*y+u3*z+u5*x2+u6*xy+u7*xz+u9*y2+u10*yz+u12*z2;
poly Hy=v1*x+v2*y+v3*z+v5*x2+v6*xy+v7*xz+v9*y2+v10*yz+v12*z2;
poly Hz=w1*x+w2*y+w3*z+w5*x2+w6*xy+w7*xz+w9*y2+w10*yz+w12*z2;
poly X=jet(psi^5
       -subst(Hx,x,(a1*t)^5,y,(a1*t)^7+a*(a1*t)^9+b*(a1*t)^11,
       z,(a1*t)^8+(a1*t)^9),13);
poly Y=jet(psi^7+m*psi^9+n*psi^11
      -subst(Hy,x,(a1*t)^5,y,(a1*t)^7+a*(a1*t)^9+b*(a1*t)^11,
      z,(a1*t)^8+(a1*t)^9),13);
poly Z=jet(psi^8+psi^9
     -subst(Hz,x,(a1*t)^5,y,(a1*t)^7+a*(a1*t)^9+b*(a1*t)^11,
     z,(a1*t)^8+(a1*t)^9),13);
matrix M1=coef(X,t);matrix M2=coef(Y,t);matrix M3=coef(Z,t);
ideal I;int ii;
for(ii=1;ii<=ncols(M1);ii++){I[size(I)+1]=M1[2,ii];}
for(ii=1;ii<=ncols(M2);ii++){I[size(I)+1]=M2[2,ii];}
for(ii=1;ii<=ncols(M3);ii++){I[size(I)+1]=M3[2,ii];}

ring S=0,(a1,a2,a3,a4,a5,a6,a7,a8,a9,a10,u1,u2,u3,u4,u5,u6,u7,u8,u9,
u10,u11,u12,u13,u14,v1,v2,v3,v4,v5,v6,v7,v8,v9,v10,v11,v12,v13,v14,
w1,w2,w3,w4,w5,w6,w7,w8,w9,w10,w11,w12,w13,w14,z1,z2,z3,z4,z5,z6,z7,
z8,z9,z10,z11,z12,z13,z14,a,b,c,d,e,f,g,h,i,j,k,l,m,n),lp;
ideal I=imap(R,I);
ideal J=std(I);
J;
> J[1]=64*a^2-63*a-56*b-64*m^2+63*m+56*n
\end{verbatim}
This is the only relation among $a,b,m,n$, which follows by inspection of the other generators of $J$, which were calculated by the elimination ordering lp. 
\end{enumerate}

\section{Proof of the Main Theorem} \label{sec6}

The proof of Theorem \ref{101} is now an easy consequence of the preceeding chapters.

\begin{proof}
Let $f=(x(t),y(t),z(t))$ define an unimodal parametrized curve singularity and $\Gamma$ its semigroup. 
Corollary \ref{cor.cho}  implies that the order of $f$ is at least $4$ and $\langle 4,7\rangle < \Gamma$. Lemma \ref{1.3.1} implies that the order of $f$ is at most $5$. If the order of $f$ is $4$ we obtain from Lemma \ref{1.3.1}  that
$\Gamma <\langle 4,9,14\rangle$ and in case of order $5$ that $\Gamma <\langle 5,7,16\rangle$. \\
Normal forms of parametrized curve singularities up to $\mathcal A$-equivalence with semigroup $\Gamma$ and
$\langle 4,7\rangle <\Gamma <\langle 4,9,14\rangle$  are given in Chapter \ref{sec4}.\\
Normal forms of parametrized curve singularities of order $5$ up to $\mathcal A$-equivalence with semigroup $\Gamma$ and 
$\Gamma <\langle 5,7,9\rangle$ are 
 given in Chapter \ref{sec5}. 
 
Since the normal forms depend only on one parameter (its modulus), they are at most unimodal. The computations also show that the modulus is free and not restricted by any relation. 
This implies that the singularities classified above are in fact unimodal.
\end{proof}

{\bf Acknowledgments:} The second author is thankful to Sultan Qaboos University for providing necessary facilities to complete this article.
We would like to
thank the Mathematical Institute of the University of Freiburg for its hospitality and the
RPTU Kaiserslautern for its support. 

This article contains no AI-generated text or mathematics.

\medskip

\smallskip\noindent
Received:  \\
Revised: \\
Accepted:

\salt

\adresa{$^{(2)}$Departement of Mathematics, College od Science, Sultan Qaboos University, \\ PO Box 36, AI Khod 123, Sultanate of Oman \\
E-mail: {\tt K.mehmood@squ.edu.om}
}
\adresa{$^{(1), (3)} $Department of Mathematics, RPTU Kaiserslautern, \\Erwin-Schr\"odinger-Str., 67663 Kaiserslautern, Germany
  \\
E-mail (1): {\tt greuel@mathematik.uni-kl.de}\\ 
E-mail (3): {\tt pfister@mathematik.uni-kl.de}
}
\end{document}